\documentclass[reqno,a4paper]{amsart}
\usepackage[T1]{fontenc}
\usepackage[utf8]{inputenc}
\usepackage[english]{babel}
\usepackage{microtype}

\usepackage{latexsym, amssymb, amsmath}
\usepackage{enumitem}
\usepackage{amsthm}
\usepackage{tikz-cd}
\usepackage{orcidlink}
\usepackage{doi}
\renewcommand{\doi}[1]{DOI:~\href{https://doi.org/#1}{#1}}
\usepackage{url}

\usepackage{verbatim}
\usepackage{calrsfs}
\usepackage{mathabx}
\usepackage{xcolor}
\usepackage{mathtools}
\usepackage{blkarray}
\usepackage[textsize=scriptsize]{todonotes}
\usepackage{thmtools}
\usepackage{xurl}
\usepackage{nameref,hyperref,cleveref}
\hypersetup{hidelinks}
\usepackage{quiver}

\theoremstyle{plain}
\newtheorem{theorem}[subsection]{Theorem}
\newtheorem{proposition}[subsection]{Proposition}
\newtheorem{corollary}[subsection]{Corollary}
\newtheorem{lemma}[subsection]{Lemma}

\theoremstyle{definition}
\newtheorem{definition}[subsection]{Definition}

\theoremstyle{remark}
\newtheorem{remark}[subsection]{Remark}
\newtheorem{example}[subsection]{Example}
\newtheorem{examples}[subsection]{Examples}

\crefname{examples}{Examples}{Examples}
\Crefname{examples}{Examples}{Examples}
\crefname{problem}{Open problem}{Open problems}
\Crefname{problem}{Open problem}{Open problems}

\numberwithin{equation}{section}

\newenvironment{tfae}
{
\begin{enumerate}}
{\end{enumerate}}

\newcommand{\noproof}{\hfil\qed}

\newcommand{\C}{\ensuremath{\mathcal{C}}}

\newcommand{\F}{\mathbb{F}}         
\newcommand{\E}{\mathcal{E}}

\newcommand{\V}{\mathcal{V}}
\newcommand{\Z}{\mathbb{Z}}
\newcommand{\U}{\mathcal{U}}

\DeclareMathOperator{\bchar}{char}
\DeclareMathOperator{\Ker}{Ker}
\DeclareMathOperator{\Aut}{Aut}
\DeclareMathOperator{\op}{op}
\DeclareMathOperator{\Act}{Act}
\DeclareMathOperator{\Hom}{Hom}
\DeclareMathOperator{\Bim}{Bim}
\DeclareMathOperator{\SplExt}{SplExt}
\DeclareMathOperator{\Imm}{Im}
\DeclareMathOperator{\End}{End}
\DeclareMathOperator{\Der}{Der}
\DeclareMathOperator{\ADer}{ADer}
\DeclareMathOperator{\Bider}{Bider}
\DeclareMathOperator{\Inn}{Inn}
\DeclareMathOperator{\M}{M}
\DeclareMathOperator{\bZ}{Z}
\DeclareMathOperator{\USGA}{USGA}
\DeclareMathOperator{\Ann}{Ann}

\DeclareMathOperator{\id}{id}

\newcommand{\Set}{\ensuremath{\mathsf{Set}}}
\newcommand{\Rng}{\ensuremath{\mathsf{Rng}}}

\newcommand{\Assoc}{\ensuremath{\mathsf{Assoc}}}

\newcommand{\Lie}{\ensuremath{\mathsf{Lie}}}
\newcommand{\Leib}{\ensuremath{\mathsf{Leib}}}
\newcommand{\Pois}{\ensuremath{\mathsf{Pois}}}
\newcommand{\CPois}{\ensuremath{\mathsf{CPois}}}
\newcommand{\BoolRng}{\ensuremath{\mathsf{BoolRng}}}
\newcommand{\BoolAlg}{\ensuremath{\mathsf{BoolAlg}}}
\newcommand{\CAssoc}{\ensuremath{\mathsf{CAssoc}}}

\newcommand{\Com}{\ensuremath{\mathsf{Com}}}
\newcommand{\ACom}{\ensuremath{\mathsf{ACom}}}
\newcommand{\JJord}{\ensuremath{\mathsf{JJord}}}
\newcommand{\Jord}{\ensuremath{\mathsf{Jord}}}
\newcommand{\AbAlg}{\ensuremath{\mathsf{AbAlg}}}
\newcommand{\Alg}{\ensuremath{\mathsf{Alg}}}
\newcommand{\PAlg}{\ensuremath{\mathsf{PAlg}}}
\newcommand{\Nil}{\ensuremath{\mathsf{Nil}}}

\newcommand{\Alt}{\ensuremath{\mathsf{Alt}}}
\newcommand{\WAlt}{\ensuremath{\mathsf{WAlt}}}
\newcommand{\Grp}{\ensuremath{\mathsf{Grp}}}

\newcommand{\Ring}{\ensuremath{\mathsf{Ring}}}
\newcommand{\CRng}{\ensuremath{\mathsf{CRng}}}
\newcommand{\CRing}{\ensuremath{\mathsf{CRing}}}

\newcommand{\SplExtX}{\ensuremath{\mathsf{SplExt}(X)}}

\date{}

\newsavebox{\pullback}
\sbox\pullback{%
\begin{tikzpicture}
\draw (0,0) -- (1ex,0ex);
\draw (1ex,0ex) -- (1ex,1ex);
\end{tikzpicture}}

\begin{document}

\title[On the representability of actions of unital algebras]{On the representability of actions\\of unital algebras}

\author[M.~Mancini]{Manuel Mancini~\orcidlink{0000-0003-2142-6193}}
\author[F.~Piazza]{Federica Piazza~\orcidlink{0009-0001-1028-9659}}
\author[C.~Vienne]{Corentin Vienne~\orcidlink{0000-0002-0397-3495}}

\email{manuel.mancini@uclouvain.be; manuel.mancini@unipa.it}
\email{federica.piazza07@unipa.it; federica.piazza1@studenti.unime.it}
\email{corentin.vienne@proton.me}

\address[M.~Mancini]{Institut de Recherche en Mathématique et Physique, Université catholique de Louvain, chemin du cyclotron 2 bte L7.01.02, B--1348 Louvain-la-Neuve, Belgium.}
\address[M.~Mancini, F.~Piazza]{Dipartimento di Matematica e Informatica, Università degli Studi di Palermo, Via Archirafi 34, 90123 Palermo, Italy.}
\address[F.~Piazza]{Dipartimento di Scienze Matematiche e Informatiche, Scienze Fisiche e Scienze della Terra, Università degli Studi di Messina, Viale Ferdinando Stagno d'Alcontres 31, 98166 Messina, Italy.}


\begin{abstract}
Working in the setting of \emph{ideally exact categories}, we investigate the representability of actions of unital non-associative algebras over a field. We show that, in general, such categories fail to be \emph{action representable}: for instance, the category of all unital algebras is not even \emph{action accessible}. We then consider this problem in the context of operadic, action accessible, unit-closed varieties. Using the construction of the \emph{external weak actor}, we prove that for any algebra $X$ in such a variety $\V$, the canonical map into its external weak actor is an isomorphism if and only if $X$ is unital. Consequently, the ideally exact category $\V_1$ of unital algebras in $\V$ is action representable, and the actor of a unital algebra $X$ is $X$ itself. The action accessibility of~$\V$ is not necessary for this conclusion: in characteristic zero, unital Jordan algebras form an action representable category, although the variety of \emph{all} Jordan algebras is not action accessible. Finally, we prove action representability for the categories of unital Poisson algebras and unital commutative Poisson algebras.
\end{abstract}

\subjclass[2020]{08A35; 08C05; 16B50; 16W25; 17A36; 17B63; 17C10; 18E13}
\keywords{Action representable category, ideally exact category, internal action, split extension, non-associative algebra, unital algebra, Jordan algebra, Poisson algebra}

\maketitle

\section{Introduction}\label{SecIntro}

The concept of \emph{internal object action} was introduced by F.~Borceux, G.~Janelidze, and G.~M.~Kelly in~\cite{BJK} with the goal of extending the correspondence between actions and split extensions from the setting of groups to the more general setting of semi-abelian categories~\cite{Semi-Ab}. In some of those categories, internal actions are exceptionally well behaved, in the sense that the actions on each object~$X$ are \emph{representable}: this means that there exists an object~$[X]$, called the \emph{actor} of $X$, such that the functor $\Act(-,X) \cong \SplExt(-,X)$ which sends an object $B$ to the set of internal actions/isomorphism classes of split extensions of $B$ on/by $X$, is naturally isomorphic to the functor $\Hom(-,[X])$. This is the case for the category of groups, where $[X]=\Aut(X)$ is the group of automorphisms of $X$. 

The study of action representability in semi-abelian categories was further developed in~\cite{BJK2}, where it was shown, for example, that the variety of commutative associative algebras over a field fails to be action representable. Later, the article~\cite{Tim} established that amongst varieties of non-associative algebras over an infinite field of characteristic different from $2$, only the category of abelian algebras and that of Lie algebras satisfy this property. The restrictive nature of action representability naturally led to the introduction of weaker, yet related, notions.

In~\cite{act_accessible}, D.~Bourn and G.~Janelidze introduced the notion of \emph{action accessibility} to encompass significant examples of categories that do not satisfy action representability, such as (not necessarily unital) rings, associative algebras, and Leibniz algebras~\cite{loday}. A.~Montoli later proved in~\cite{Montoli} that every \emph{Orzech category of interest}~\cite{Orzech} is action accessible. Additionally, the authors of~\cite{Casas} presented a broader notion of representing object in any Orzech category of interest: the \emph{universal strict general actor}.

More recently, G.~Janelidze introduced in~\cite{WRA} the notion of \emph{weak action representability}. This condition requires for each object $X$ of a semi-abelian category the existence of a \emph{weak representation} of the actions on $X$, i.e., an object $T(X)$ along with a natural monomorphism of functors $\SplExt(-,X) \rightarrowtail \Hom(-,T(X))$. Examples of weakly action representable categories are the variety of associative algebras~\cite{WRA} and the variety of Leibniz algebras~\cite{CigoliManciniMetere}.

In~\cite{WRAAlg}, the authors explored the concept of weak action representability within the framework of varieties of non-associative algebras over a field. They worked towards the construction of an external weakly representing object $\E(X)$ for actions on an object $X$ of a variety of non-associative algebras $\V$. They actually obtained a \emph{partial algebra} $\E(X)$, called the \emph{external weak actor} of $X$, together with a monomorphism of functors $\SplExt(-,X) \rightarrowtail \Hom_\PAlg(\widetilde{U}(-),\E(X))$, where $\PAlg$ is the category of partial algebras and $\widetilde{U} \colon \V \to \PAlg$ denotes the forgetful functor. In the present article, we adopt the formulation of partial algebras of~\cite{SurveyActions}, which modifies the one of~\cite{WRAAlg}.

G.~Janelidze later extended the notions of action accessibility and (weak) action representability to the broader setting of \emph{ideally exact categories}~\cite{IdeallyExact}, which were introduced as a generalisation of semi-abelian categories in order to include relevant examples of non-pointed categories, such as the categories $\Ring$ and $\CRing$ of unital rings and unital commutative rings, all non-trivial protomodular varieties of universal algebras, and any cotopos. In this context, the action representability of an ideally exact category~$\U$ is defined relative to the semi-abelian category $(\U \downarrow 0)$, where~$0$ denotes the initial object of $\U$; see \Cref{def_act_repr}. It was also shown in~\cite{IdeallyExact} that both $\Ring$ and $\CRing$ are action representable in this broader sense.

The aim of this manuscript is to study the representability of actions within the setting of categories of unital non-associative algebras. After recalling the necessary background in Sections \ref{SecVar} and \ref{SecAct}, we show in \Cref{SecIdEx} that in general a category of unital algebras fails to be action representable. In fact, we prove that the categories of unital algebras and unital commutative algebras (both not necessarily associative) are not even action accessible. 

We then study in \Cref{SecMain} the representability of actions of unital algebras in any operadic, action accessible, unit-closed variety $\V$. Using the construction of the external weak actor $\E(X)$ of an algebra $X$ of $\V$, we prove that the canonical partial algebra homomorphism $\Inn \colon X \to \E(X)$, defined by $\Inn(x) = (L_x,R_x)$, where $L_x$ and $R_x$ denote the left and right multiplications by $x$, is an isomorphism if and only if~$X$ is unital. In \Cref{final_thm}, the main result of this article, we prove that the ideally exact category $\V_1$ of unital algebras in $\V$ is action representable, and the actor of a unital algebra $X$ is $X$ itself. This result implies that, in particular, the categories of unital alternative, unital associative, and unital commutative associative algebras are action representable, and highlights a contrast with the non-unital case, where action representability is a highly restrictive property, essentially characterising the variety of Lie algebras~\cite{Tim}. We also show in \Cref{prop_jord} that the action accessibility of $\V$ is not necessary for $\V_1$ to be action representable: over a field of characteristic zero, the variety of Jordan algebras is unit-closed and not action accessible, while the ideally exact category of unital Jordan algebras is action representable, the actor of a unital Jordan algebra $X$ being $X$ itself. To the best of our knowledge, this is the first example of an action representable ideally exact category whose associated semi-abelian category is not even action accessible.

Finally, in \Cref{SecPois}, over a field of characteristic different from~$2$, we study the representability of actions of the categories $\Pois_1$ and $\CPois_1$ of unital Poisson algebras and unital commutative Poisson algebras. Here we employ the explicit description of split extensions and derived actions of Poisson algebras given in~\cite{CigoliManciniMetere}, and the construction of the external weak actor given in~\cite{SeveralMult}. While the categories $\Pois$ and $\CPois$ of Poisson algebras and commutative Poisson algebras are not action representable~\cite{SeveralMult}, and it remains open whether they have weakly representable actions, we prove that the ideally exact categories $\Pois_1$ and $\CPois_1$ are action representable, and that the actor of a unital (commutative) Poisson algebra $X$ is isomorphic to $X$ itself.

\section{Varieties of non-associative algebras}\label{SecVar}

We now describe the algebraic framework in which we work: \emph{varieties of non-associative algebras} over a field $\F$. We think of those as collections of algebras satisfying a chosen set of polynomial equations. We refer the reader to~\cite{VdL-NAA} for more details.

A \emph{non-associative algebra} over $\F$ is a vector space $X$ equipped with a bilinear multiplication $X \times X \to  X$, denoted by $(x,y) \mapsto xy$. In general, the existence of a multiplicative identity is not assumed. The category of \emph{all} non-associative algebras over $\F$ is denoted by $\Alg$ and it has as morphisms the linear maps that preserve the multiplication.

\begin{definition}\label{def identity variety}
An \emph{identity} of an algebra $X$ is a non-associative polynomial
\[
\varphi=\varphi(x_1,\ldots,x_n)
\]
such that $\varphi(x_1,\ldots,x_n)=0$ for any $x_1,\ldots,x_n \in X$. We say that the algebra $X$ \emph{satisfies} the identity $\varphi$.
\end{definition}

If $I$ is a set of identities, then the \emph{variety of non-associative algebras} $\V$ determined by $I$ is the class of all algebras which satisfy all the identities in $I$. Conversely, we say that a variety of non-associative algebras $\V$ satisfies the identities in $I$ if every algebra of $\V$ satisfies the given set of identities. 

\begin{remark}
From now on, when we consider a variety $\V$, we assume that $\V$ is a variety of non-associative algebras over a field $\F$. We fix the field $\F$, so that we may drop it from our notation.
\end{remark}

For a monomial~$\varphi=\varphi(x_1,\ldots,x_n)$, we define its \emph{type} as the $n$-tuple $(k_1, \ldots , k_n)$, where $k_i$ is the number of times $x_i$ appears in $\varphi$, and its \emph{degree} as the natural number $k_1 + \cdots + k_n$. A polynomial is said to be \emph{homogeneous} if all its monomials are of the same type, and it is said to be \emph{multilinear} if all its monomials have type~$(1,\ldots,1)$.

\begin{definition}
A variety of non-associative algebras $\V$ is said to be \emph{homogeneous} (respectively \emph{operadic}) if it is determined by a set of homogeneous (respectively multilinear) identities.
\end{definition}

We recall that, if the field $\F$ is infinite, then any variety of non-associative $\F$-algebras $\V$ is homogeneous~\cite{hmg}. Furthermore, if $\bchar(\F)=0$, then any variety $\V$ is operadic. This follows from the multilinearisation process, see~\cite[Corollary 3.7]{Osborn}.

\begin{examples}\label{Examples varieties}{\ }
\begin{enumerate}
\item $\AbAlg$ is the variety of \emph{abelian} algebras, which is determined by the identity~$xy=0$.

\item $\Com$ and $\ACom$ denote, respectively, the varieties of commutative and anti-commutative (not necessarily associative) algebras. These varieties are determined, respectively, by the identities $xy-yx=0$ and $xy+yx=0$.

\item $\Assoc$ is the variety of \emph{associative} algebras, which is determined by the identity $x(yz)-(xy)z=0$.

\item $\CAssoc$ denotes the variety of commutative associative algebras.

\item $\Alt$ is the variety of \emph{alternative algebras}, which is determined by the identities
\[
(yx)x-y(xx)=0 \quad \text{and} \quad x(xy)-(xx)y=0.
\]
If $\bchar(\F) \neq 2$, then, as described in Example~\ref{example Alt}, the variety $\Alt$ can be interpreted as an operadic variety.

We observe that every associative algebra is alternative, while an example of an alternative algebra which is not associative is given by the algebra of octonions.

\item $\Lie$ denotes the variety of \emph{Lie algebras}, characterised by $x^2=0$ together with the \emph{Jacobi identity} $x(yz)+y(zx)+z(xy)=0$. We recall that when $\bchar(\F) \neq 2$, the condition $x^2=0$ is equivalent to anti-commutativity, and hence in this situation the variety $\Lie$ is operadic.

\item $\Jord$ is the variety of \emph{Jordan algebras}~\cite{Jacobson}, which is determined by commutativity and the \emph{Jordan identity} $x^2(xy)=x(x^2y)$. If $\bchar(\F)\neq2$, every associative algebra becomes a Jordan algebra with respect to the product $a\circ b=(ab+ba)/2$.

\item $\JJord$ is the variety of \emph{Jacobi--Jordan algebras}~\cite{Burde}, which is determined by commutativity and the Jacobi identity. Jacobi--Jordan algebras, also known as \emph{mock-Lie algebras}, are the commutative counterpart of Lie algebras. Every Jacobi--Jordan algebra is a Jordan algebra, which justifies the name. Indeed, the Jacobi identity applied to $x$, $x$ and $y$ gives $x^2y=-2x(xy)$, so that $x^2(xy)=-2x(x(xy))=x(x^2y)$.

\item $\Leib$ is the variety of \emph{(right) Leibniz algebras}, which is determined by the \emph{(right) Leibniz identity}, that is $(xy)z-(xz)y-x(yz)=0$.

\item $\BoolRng$ is the category of (not necessarily unital) \emph{Boolean rings}, which may be seen as the variety of commutative associative $\F_2$-algebras satisfying $x^2=x$. We remark that this variety is not homogeneous.

\item Given any variety $\V$, one can define the subvariety $\Nil_k(\V)$ of $k$-\emph{nilpotent} algebras of $\V$, which is determined by all the identities of $\V$ together with
\[
x_1 x_2 \cdots x_{k+1}=0,
\]
for all possible choices of parentheses. For instance, if $\V=\CAssoc$, every non-trivial multilinear identity of degree $k+1$ is equivalent to $x_1\cdots x_{k+1}=0$. Hence, the only proper operadic subvarieties of $\CAssoc$ are the varieties $\Nil_k(\CAssoc)$ of $k$-nilpotent commutative associative algebras, for $k\geq1$, where $\Nil_1(\CAssoc)=\AbAlg$, together with the trivial variety, which corresponds to $k=0$.
\end{enumerate}
\end{examples}

In this paper, we aim to study the representability of actions of unital algebras. To do this, we need to recall the following definition.

\begin{definition}\cite{ShestakovBittencourt}
A variety $\V$ is said to be \emph{unit-closed} if for any algebra $X$ of~$\V$, the algebra $\widetilde{X}=\langle X,1 \rangle$ obtained by adjoining to $X$ an external unit $1$ which satisfies the identities $x \cdot 1 = 1 \cdot x = x$ for any $x \in X$, is still an object of~$\V$. 
\end{definition}

For instance, the varieties $\Com$, $\Assoc$, $\CAssoc$, $\Alt$, $\Jord$, and $\BoolRng$ are unit-closed, whereas the varieties $\AbAlg$, $\Leib$, and, more generally, any variety of anti-commutative algebras over a field of characteristic different from $2$, are not unit-closed. Note that in characteristic $2$, commutativity and anti-commutativity coincide. Thus, the property of being unit-closed depends on the identities satisfied by the variety, and the characteristic of the base field plays a crucial role. For example, the variety $\JJord$ is unit-closed if and only if $\bchar(\F)=3$.

Let us check that $\Jord$ is unit-closed, the other cases being straightforward. For an element $a$ of an algebra $A$, we denote by $L_a\colon A\to A$ the left multiplication by~$a$ and, for linear endomorphisms $f$ and $g$ of a vector space, we write $[f,g]=f\circ g-g\circ f$. Thus, a commutative algebra $A$ is a Jordan algebra if and only if $[L_{x^2},L_x]=0$ for every $x\in A$. Let $X$ be a Jordan algebra, and let $x=\lambda 1+a$ be an element of $\widetilde X$, with $\lambda\in\F$ and $a\in X$. Then $L_x=\lambda\id+L_a$ and $L_{x^2}=\lambda^2\id+2\lambda L_a+L_{a^2}$, so that $[L_{x^2},L_x]=[L_{a^2},L_a]$. This commutator vanishes on $X$, since $X$ is a Jordan algebra, and on $1$, since $a^2a=aa^2$ by commutativity. Since $\widetilde X$ is commutative, it is a Jordan algebra.

\begin{remark}
Given a unit-closed variety $\V$, we denote by $\V_1$ the subclass of unital algebras of $\V$. We will see in \Cref{SecIdEx} that $\V_1$ is an \emph{ideally exact category} in the sense of G.~Janelidze~\cite{IdeallyExact}.
\end{remark}

\section{Representability of actions in semi-abelian categories}\label{SecAct}

Semi-abelian categories were introduced in~\cite{Semi-Ab} by G.~Janelidze, L.~Márki and W.~Tholen in order to provide a categorical setting which would capture categorical-algebraic properties of groups, rings and algebras. We recall that a category~$\C$ is semi-abelian when it is pointed, Barr-exact~\cite{Barr}, Bourn-protomodular~\cite{Bourn} and has finite coproducts. We refer to the monograph~\cite{malcev} for a detailed account of these notions.

\begin{example}
Any variety of non-associative algebras $\V$ can be seen as a category where the morphisms are the same ones as in $\Alg$ and we have a full inclusion functor $\V \hookrightarrow \Alg$. In particular, any such variety is a semi-abelian category, see~\cite[Theorem~9.5]{VdL-NAA}.
\end{example}

A notion that can be explored in the context of semi-abelian categories is that of split extension.

\begin{definition}
A \emph{split extension} of an object $B$ by another object $X$ of a semi-abelian category $\C$ is a diagram
\begin{equation*}\label{eq:split_ext}
\begin{tikzcd}
X\arrow [r, "k"]
&A \arrow[r, shift left, "\alpha"] &
B\ar[l, shift left, "\beta"]
\end{tikzcd}
\end{equation*}
in $\C$, such that $\alpha \circ \beta = \id_B$ and $(X,k)$ is a kernel of $\alpha$.
\end{definition}

For any object $X$ of $\C$, one can consider the functor
\[
\SplExt(-,X) \colon \C^{\op} \to \Set
\]
which assigns to each object $B$ of $\C$ the set $\SplExt(B,X)$ of isomorphism classes of split extensions of $B$ by $X$, and to any morphism $f \colon B'\to  B$ the \emph{change of base} map $f^* \colon \SplExt(B,X) \to \SplExt(B',X)$ given by pulling back along $f$.

In a semi-abelian category $\C$, one can also define the notion of internal action~\cite{BJK}. 
Internal actions on an object $X$ give rise to a functor
\[
\Act(-,X)\colon \C^{\op} \to \Set,
\]
and it was proved in~\cite{BJK} that there is a natural isomorphism
\[
\Act(-,X)\cong \SplExt(-,X).
\]
We do not describe internal actions in detail here, since split extensions provide a more practical framework, especially in the context of non-associative algebras.

\begin{definition}\cite{BJK}
A semi-abelian category $\C$ is said to be \emph{action representable} if the functor $\SplExt(-,X)$ is representable for every object $X$ of $\C$. That is, for every object $X$ there exists an object $[X]$ of $\C$ along with a natural isomorphism of functors
\[
\SplExt(-,X) \cong \Hom_{\C}(-,[X]).
\]
\end{definition}

Examples of action representable categories are the category $\Grp$ of groups and the category $\Lie$ of Lie algebras over a unital commutative ring. In the case of groups, every action of $B$ on $X$ corresponds to a group homomorphism $B \to \Aut(X)$, where $\Aut(X)$ is the automorphism group of $X$. Similarly, for Lie algebras, every action of~$B$ on $X$ is described by a Lie algebra homomorphism $B \to  \Der(X)$, where $\Der(X)$ is the Lie algebra of derivations of $X$.

However, action representability is a rather restrictive property. For instance, it was shown in~\cite[Theorem~6.1]{Tim} that amongst varieties of non-associative algebras over an infinite field $\F$ with $\bchar(\F) \neq 2$, the only action representable varieties are the category of abelian algebras and that of Lie algebras.

\begin{remark}\label{BJ}
It is shown in~\cite{BJK} that action representability is equivalent to the condition that, for every object $X$ in $\C$, the category $\SplExtX$ of split extensions with kernel $X$ has a terminal object of the form
\[
\begin{tikzcd}
X \arrow[r]
& {[X] \ltimes X} \arrow[r, shift left] &
{[X]} \arrow[l, shift left]
\end{tikzcd}
\]
\end{remark}

We can weaken this condition by assuming instead that every object in $\SplExtX$ is \emph{accessible}, i.e., it has a morphism into a \emph{subterminal} or \emph{faithful} object~\cite{act_accessible}, that is, an object which admits at most one morphism into it. This leads to the following notion, introduced by D.~Bourn and G.~Janelidze in~\cite{act_accessible} in order to compute centralisers of normal subobjects or of equivalence relations.

\begin{definition}\cite{act_accessible}
A semi-abelian category $\C$ is \emph{action accessible} if, for any object~$X$ of $\C$, every object in $\SplExtX$ is accessible.
\end{definition}

The notion of action accessibility allows us to encompass a wider class of categories that are not action representable, such as the categories $\Rng$ and $\CRng$ of not necessarily unital rings and not necessarily unital commutative rings, the varieties $\Assoc$ and $\CAssoc$, and the variety $\Pois$ of Poisson algebras.

\begin{remark}
Since the existence of a terminal object in $\SplExtX$ is stronger than every object being accessible, it follows from \Cref{BJ} that
\[
\emph{action representability}\Rightarrow\emph{action accessibility.}
\]
\end{remark}

Alternatively, instead of weakening the condition on the existence of a terminal object in $\SplExtX$, one may weaken the properties of the representing object $[X]$. This is the approach taken in~\cite{Casas}, where the authors showed that every \emph{Orzech category of interest}~\cite{Orzech} admits a \emph{universal strict general actor} (USGA for short). A.~Montoli later proved in~\cite{Montoli} that all Orzech categories of interest are action accessible.

Recently in~\cite{WRA}, G.~Janelidze introduced the notion of \emph{weakly representable action}.

\begin{definition}\cite{WRA}
A semi-abelian category $\C$ is said to be \emph{weakly action representable} if the functor $\SplExt(-,X)$ admits a weak representation for every object~$X$ of~$\C$. This means that for every object $X$, there exists an object $T=T(X)$ of $\C$ and a natural monomorphism of functors
\[
\mu \colon \SplExt(-,X) \rightarrowtail \Hom_{\C}(-,T).
\]
A morphism $\varphi\colon B \to  T$ that belongs to $\Imm(\mu_B)$ is called an \emph{acting morphism}.
\end{definition}

Examples of weakly action representable categories include the variety $\Assoc$ of associative algebras (see~\cite{WRA}), where $T=\Bim(X)$ is the associative algebra of \emph{bimultipliers} of $X$~\cite{MacLane58}; the variety $\Leib$ of Leibniz algebras (see~\cite[Theorem~3.14]{CigoliManciniMetere}), where $T=\Bider(X)$ is the Leibniz algebra of \emph{biderivations} of $X$~\cite{loday}; and the varieties $\Nil_2(\Com)$ and $\Nil_2(\ACom)$ of $2$-nilpotent commutative and anti-commutative algebras (see Theorems 2.19 and 2.21 of~\cite{WRAAlg}). Moreover, if in an Orzech category of interest $\C$ each $\USGA(X)$ is an object of the category, then $\C$ is weakly action representable (see~\cite[Corollary 4.2]{CigoliManciniMetere}).

A result obtained by G.~Janelidze in~\cite{WRA} is that every weakly action representable category is action accessible. Thus, we have
\[
\textit{action representability}\Rightarrow\textit{weak action representability}\Rightarrow\textit{action accessibility.}
\]
The converse of the first implication does not hold, since the variety $\Assoc$ of associative algebras is not action representable (see~\cite[Proposition 1.11]{Tim}). J.~R.~A.~Gray recently observed in~\cite{Gray} that the converse of the second implication also fails: he proved that the varieties of $k$-solvable groups ($k \geq 3$) are not weakly action representable but action accessible. This result was then extended in~\cite{XabiMancini} to the varieties of $2$-solvable groups and $n$-nilpotent groups ($n \geq 3$). There, the authors also proved that the varieties of $k$-solvable Lie algebras ($k \geq 2$) and $n$-nilpotent Lie algebras ($n \geq 3$) fail to be weakly action representable.

\subsection{External weak actor in varieties of non-associative algebras}

The representability of actions of non-associative algebras was studied in~\cite{WRAAlg}, where the authors proved that, for any object $X$ of an operadic, action accessible variety $\V$, there exists a \emph{partial algebra} $\E(X) \subseteq \End(X) \times \End(X)$, called the \emph{external weak actor} of $X$, together with a monomorphism of functors
\[
\tau \colon \SplExt(-,X) \rightarrowtail \Hom_\PAlg(\widetilde{U}(-),\E(X)),
\]
called the \emph{external weak representation}. Here, $\PAlg$ denotes the category of \emph{partial algebras}, and $\widetilde{U} \colon \V \to \PAlg$ is the forgetful functor.

We now recall the notion of partial algebra, following the formulation of~\cite{SurveyActions}, which modifies the one of~\cite{WRAAlg}, where the domain $\Omega\subseteq X\times X$ of a partial multiplication on a vector space $X$ was required to be a vector subspace. Since products with zero are always defined in the constructions under consideration, this forces $\Omega=X\times X$, so that the multiplication is necessarily total. To allow non-total multiplications, the domain is instead taken to be a vector subspace of $X\otimes X$, as follows.

\begin{definition}\cite{SurveyActions}
A \emph{partial algebra} over $\F$ is a triple $(X,D,\mu)$, where $X$ is an $\F$-vector space, $D\leq X\otimes X$ is a vector subspace and $\mu\colon D\to X$ is a linear map. Given $x,y\in X$, we say that the product $xy$ is \emph{defined} if $x\otimes y\in D$, and in this case we set $xy\coloneqq\mu(x\otimes y)$. The partial algebra $(X,D,\mu)$ is said to be \emph{total} if $D=X\otimes X$.

A \emph{homomorphism of partial algebras} $(X,D,\mu)\to(X',D',\mu')$ is a linear map $f\colon X\to X'$ such that
\[
(f\otimes f)(D)\subseteq D' \quad \text{and} \quad f\circ\mu=\mu'\circ\left.(f\otimes f)\right|_D.
\]
\end{definition}

We denote by $\PAlg$ the category of partial algebras over $\F$ and their homomorphisms. Every non-associative algebra $X$ determines the total algebra
\[
\widetilde U(X)=(X,X\otimes X,\widetilde m_X),
\]
where $\widetilde m_X\colon X\otimes X\to X$ is the linearisation of the multiplication of $X$. This defines a fully faithful functor $\widetilde U\colon\V\to\PAlg$ for every variety $\V$. Moreover, if~$B$ is a non-associative algebra and $(X',D',\mu')$ is a partial algebra, a linear map $\varphi\colon B\to X'$ defines a partial algebra homomorphism $\widetilde U(B)\to(X',D',\mu')$ if and only if $\varphi(b)\otimes\varphi(b')\in D'$ and $\varphi(bb')=\mu'(\varphi(b)\otimes\varphi(b'))$ for all $b,b'\in B$ (see~\cite{SurveyActions}). Finally, an element $e$ of a partial algebra $(X,D,\mu)$ is called a \emph{unit} if $e\otimes x,\,x\otimes e\in D$ and $ex=xe=x$ for every $x\in X$. As in the case of total algebras, a unit is unique when it exists: if $e$ and $e'$ are units, then $e\otimes e'\in D$ and $e=ee'=e'$.

\begin{lemma}\label{lem_iso_total}
Let $B$ be a non-associative algebra and let $\varphi\colon\widetilde U(B)\to(X',D',\mu')$ be a bijective homomorphism of partial algebras. Then $D'=X'\otimes X'$, and $\varphi$ is an isomorphism of partial algebras. In particular, $\mu'$ is the linearisation of a bilinear multiplication on $X'$, and $\varphi\colon B\to X'$ is an algebra isomorphism.
\end{lemma}

\begin{proof}
Since $\varphi$ is bijective, so is $\varphi\otimes\varphi\colon B\otimes B\to X'\otimes X'$. Hence $X'\otimes X'=(\varphi\otimes\varphi)(B\otimes B)\subseteq D'$, so that $D'=X'\otimes X'$ and $\mu'=\varphi\circ\widetilde m_B\circ(\varphi\otimes\varphi)^{-1}$. It follows that $\varphi^{-1}$ is a homomorphism of partial algebras, and that $\varphi$ is an algebra isomorphism when $X'$ is endowed with the multiplication induced by~$\mu'$.
\end{proof}

We recall the construction of $\E(X)$, since it will be central in proving that actions on a unital algebra $X$ are representable. We begin by recalling a result that links action accessibility with the identities of degree $3$ satisfied by a variety $\V$.

\begin{theorem}\cite{Tim, XabiTim, Coco_Thesis}\label{Theorem AC iff Orzech}
Let $\V$ be a homogeneous variety of non-associative algebras. The following conditions are equivalent:
\begin{tfae}
\item $\V$ is an Orzech category of interest;
\item $\V$ is action accessible;
\item There exist $\lambda_{1}, \dots, \lambda_{8}, \mu_{1}, \dots, \mu_{8} \in \F$ such that
\begin{align}
\begin{split}\label{eq00}
x(yz)=\lambda_1(xy)z & +\lambda_2(yx)z+\lambda_3 z(xy) + \lambda_4 z(yx)\\
& + \lambda_5 (xz)y + \lambda_6 (zx)y + \lambda_7 y(xz) + \lambda_8 y(zx)
\end{split}
\intertext{and}
\begin{split}\label{eq01}
(yz)x=\mu_1(xy)z & +\mu_2(yx)z+\mu_3 z(xy) + \mu_4 z(yx)\\
& + \mu_5 (xz)y + \mu_6 (zx)y + \mu_7 y(xz) + \mu_8 y(zx)
\end{split}
\end{align}
are identities of $\V$. The two previous identities together are called the \emph{$\lambda/\mu$-rules}. \noproof
\end{tfae}
\end{theorem}

\begin{remark}
Note that the choice of the $\lambda/\mu$-rules is not unique, but fixed for our purposes.
\end{remark}

Now, let $\V$ be an action accessible variety determined by a set of multilinear identities
\[
\Phi_{k,i}(x_1,\ldots,x_k)=0, \quad i=1,\ldots,n_k,
\]
where $k=\deg\Phi_{k,i}\geq 2$ (an identity of degree $1$ only occurs in the trivial variety, which we exclude). For each $k$, the integer $n_k$ may be chosen finite, since the space of multilinear polynomials in $k$ fixed variables is finite-dimensional. We also fix $\lambda_1, \ldots,\lambda_8,\mu_1, \ldots,\mu_8 \in \F$ which determine a choice of the $\lambda/\mu$-rules. 

The external weak actor $\E(X)$ is defined as the subspace of all pairs $f=(f\ast-,-\ast f) \in \End(X) \times \End(X)$ satisfying
\[
\Phi_{k,i}(\alpha_1,\ldots,\alpha_k)=0, \quad i=1,\ldots,n_k,
\]
for every $k$ and every choice of $\alpha_j=f$ and $\alpha_t \in X$, with $t\neq j\in \{1,\ldots,k\}$, where we write $fx\coloneqq f\ast x$ and $xf\coloneqq x\ast f$. Since $f$ occurs exactly once, each such evaluation is an element of~$X$.

The chosen $\lambda/\mu$-rules induce a bilinear map
\[
m_X\colon\E(X)\times\E(X)\to\End(X)\times\End(X)\colon (f,g)\mapsto h=(h\ast-,-\ast h),
\]
where
\begin{align*}
x \ast h=\lambda_1(x\ast f)\ast g & +\lambda_2(f\ast x)\ast g+\lambda_3g\ast(x\ast f) + \lambda_4 g \ast(f\ast x) \\
 & + \lambda_5 (x\ast g)\ast f + \lambda_6 (g\ast x)\ast f + \lambda_7 f\ast (x\ast g) + \lambda_8 f\ast (g\ast x)
\intertext{and}
h \ast x=\mu_1(x\ast f)\ast g     & +\mu_2(f\ast x)\ast g+\mu_3g\ast(x\ast f) + \mu_4 g\ast(f\ast x) \\
 & + \mu_5 (x\ast g)\ast f + \mu_6 (g\ast x)\ast f + \mu_7 f\ast (x\ast g) + \mu_8 f\ast (g\ast x),
\end{align*}
for every $x\in X$. In general, $m_X(f,g)$ need not belong to $\E(X)$. Denoting by $\widetilde m_X\colon\E(X)\otimes\E(X)\to\End(X)\times\End(X)$ the linearisation of $m_X$, we set
\[
D_X\coloneqq\widetilde m_X^{-1}\bigl(\E(X)\bigr)\leq\E(X)\otimes\E(X)
\qquad\text{and}\qquad
\mu_X\coloneqq\left.\widetilde m_X\right|_{D_X}\colon D_X\to\E(X).
\]
Then $(\E(X),D_X,\mu_X)$ is a partial algebra, which we still denote by $\E(X)$. Its partial multiplication is denoted by $\langle-,-\rangle$: for $f,g\in\E(X)$, the product $\langle f,g\rangle$ is defined if and only if $f\otimes g\in D_X$, or equivalently $m_X(f,g)\in\E(X)$, and in this case $\langle f,g\rangle=m_X(f,g)$. Finally, $\E(X)$ is called the \emph{external actor} of $X$ whenever $\tau$ is a natural isomorphism.

Notice that $D_X$ contains the kernel of $\widetilde m_X$, which in general includes tensors that are not pure. The tensors of $D_X$ which are not linear combinations of pure tensors belonging to $D_X$ do not affect the defined products, but they do affect the homomorphisms of partial algebras. This has no consequences for the results of this article, since the homomorphisms of partial algebras used to represent actions have a total algebra as domain, and hence they only depend on the defined products of~$\E(X)$.

\begin{remark}\label{rem_one_f}
The subspace $\E(X)$ does not depend on the chosen defining family of identities. Indeed, since the multilinear identities of $\V$ are linear combinations of the $\Phi_{k,i}$ evaluated at monomials with pairwise disjoint sets of variables and multiplied by further monomials, every multilinear identity of $\V$ vanishes when exactly one of its variables is replaced by an element of $\E(X)$ and the others by elements of~$X$; see~\cite{SurveyActions}. In particular, this applies to the $\lambda/\mu$-rules.
\end{remark}

\begin{example}\label{Example associative}
If $\V=\Assoc$ and we fix the standard choice of the $\lambda/\mu$-rules (that is $\lambda_1 = \mu_8 = 1$ and $\lambda_i=\mu_j=0$, for any $i \neq 1$ and $j \neq 8$), then $\E(X)$ is a total algebra, isomorphic to the associative algebra $\Bim(X)$ of \emph{bimultipliers} of $X$ (see~\cite{MacLane58}, where they are called \emph{bimultiplications}), that is, of the pairs $(f \ast -,- \ast f) \in \End(X)\times \End(X)^{\text{op}}$ such that
\[
f \ast (xy)=(f \ast x)y,\qquad (xy) \ast f=x(y \ast f),\qquad x(f \ast y)=(x \ast f)y
\]
for all $x,y\in X$, where the multiplication is given by
\[
\langle(f\ast-,-\ast f),(g\ast-,-\ast g)\rangle
=
\bigl(f\ast(g\ast-),(-\ast f)\ast g\bigr).
\]
It was proved in~\cite{WRA} that $\Bim(X)$ provides a weak representation of the actions on~$X$.
\end{example}

\begin{example}
Let $\V=\CAssoc$, and choose the $\lambda/\mu$-rules
\[
x(yz)=(xz)y \qquad\text{and}\qquad (yz)x=(xz)y,
\]
that is, $\lambda_5=\mu_5=1$ and all the other coefficients equal to zero. Then every $f\in\E(X)$ satisfies $f\ast x=x\ast f$, and $\E(X)$ is a total algebra, which is isomorphic to the associative algebra
\[
\M(X)=\{ f \in \End(X) \mid f(xy)=f(x)y, \; \forall x,y \in X \}
\]
of \emph{multipliers} of $X$, whose multiplication is given by the usual composition of functions. Moreover, it was proved in~\cite[proof of Theorem~2.6]{BJK2} that there is a natural isomorphism
\[
\SplExt(-,X) \cong \Hom_\Assoc\left(\overline{U}(-),\M(X)\right),
\]
where $\overline{U} \colon \CAssoc \to \Assoc$ denotes the forgetful functor. Thus, $\E(X) \cong \M(X)$ is an external actor of $X$.
\end{example}

\begin{example}
Let $\bchar(\F) \neq 2$. If $\V=\Lie$, then $\E(X) \cong \Der(X)$ is the actor of~$X$.
\end{example}

\begin{remark}\label{rem_int}
Let $X$ be an algebra of $\V$, and let $L_x = x \cdot -$ and $R_x = - \cdot x$ denote, respectively, the left and right multiplications by $x\in X$. Since the defining identities of $\V$ hold in $X$, the pair $(L_x,R_x)$ belongs to $\E(X)$ for every $x\in X$. Moreover, the $\lambda/\mu$-rules, applied in $X$, give
\[
m_X\bigl((L_x,R_x),(L_y,R_y)\bigr)=(L_{xy},R_{xy})\in\E(X)
\]
for all $x,y\in X$. Hence the linear map
\[
\Inn \colon \widetilde{U}(X) \to \E(X) \colon x \mapsto (L_x,R_x)
\]
is a homomorphism of partial algebras. The kernel of $\Inn$ is the \emph{annihilator}
\[
\Ann(X)=\{ x \in X \mid xy = yx = 0,\ \forall y \in X \}
\]
of $X$. In some cases, such as for Lie algebras, $\Ann(X)$ is called the \emph{centre} of $X$ and is denoted by $\bZ(X)$.
\end{remark}

The following result explains how the construction of the external weak actor allows one to represent actions on, or equivalently split extensions by, an algebra $X$.

\begin{theorem}\cite{WRAAlg,SurveyActions}\label{thmV}
Let $\V$ be an operadic, action accessible variety of non-associative algebras, and fix a choice of the $\lambda/\mu$-rules.
\begin{enumerate}
\item Let $X$ be an object of $\V$. There exists a natural monomorphism of functors
\[
\tau \colon \SplExt(-,X) \rightarrowtail \Hom_{\PAlg}(\widetilde{U}(-),\E(X)),
\]
where $\widetilde{U} \colon \V \to \PAlg$ denotes the forgetful functor and, for every object $B$ of~$\V$, $\tau_B$ is the injection which sends a split extension
\[
\begin{tikzcd}
X\arrow [r, "k"]
&A \arrow[r, shift left, "\alpha"] &
B\ar[l, shift left, "\beta"]
\end{tikzcd}
\]
to the partial algebra homomorphism
\[
\widetilde{U}(B) \to \E(X) \colon b \mapsto (b \ast -, - \ast b),
\]
where $b \ast x \coloneqq \beta(b)k(x)$ and $x \ast b \coloneqq k(x)\beta(b)$.
\item Let $B, X$ be objects of $\V$. The partial algebra homomorphism
\[
\widetilde{U}(B) \to \E(X) \colon b \mapsto (b \ast -,- \ast b)
\]
belongs to $\Imm(\tau_B)$ if and only if $\Phi_{k,i}(\alpha_1,\ldots,\alpha_k)=0$ for every $k$ and every $i=1,\ldots,n_k$, where at least one of the $\alpha_1, \ldots, \alpha_k$ is an element $x \in X$, and the others are elements $b \in B$. Here, we write $bx\coloneqq b\ast x$ and $xb\coloneqq x\ast b$.
\item If $D_X=\E(X)\otimes\E(X)$ and the resulting total algebra, with multiplication $\langle -,- \rangle$, is an object of $\V$, then $(\E(X),\tau)$ is a weak representation of $\SplExt(-,X)$.
\noproof
\end{enumerate}
\end{theorem}

\begin{example}
Let $X$ be an object of $\V$, and consider the canonical action of $X$ on itself given by left and right multiplication, that is
\[
b \ast x = bx, \quad \text{and} \quad x \ast b = xb.
\]
Then the component $\tau_X$ sends this action to the partial algebra homomorphism $\Inn \colon \widetilde{U}(X) \to \E(X)$ of \Cref{rem_int}.
\end{example}

We conclude this section with the following proposition, showing that there are special cases where $\E(X)$ becomes the actor of $X$.

\begin{proposition}\label{prop_e(x)}
Let $\V$ be an operadic, action accessible variety of non-associative algebras, and fix a choice of $\lambda/\mu$-rules. For any object $X \in \V$, if the homomorphism of partial algebras $\Inn \colon \widetilde{U}(X) \to \E(X)$ is an isomorphism, then the actions on $X$ are representable and $X$ is its own actor.
\end{proposition}

\begin{proof}
Suppose that $\Inn$ is an isomorphism of partial algebras. Since $\widetilde U$ is fully faithful, composition with $\Inn^{-1}$ gives natural isomorphisms
\[
\Hom_{\PAlg}(\widetilde{U}(-),\E(X)) \cong \Hom_{\PAlg}(\widetilde{U}(-),\widetilde{U}(X)) \cong \Hom_{\V}(-,X),
\]
and hence, by composition with $\tau$, we get a natural monomorphism of functors
\[
\widetilde{\tau} \colon \SplExt(-,X) \rightarrowtail \Hom_{\V}(-,X).
\]
It remains to show that, for any algebra $B$ of $\V$, the component $\widetilde{\tau}_B$ is surjective, that is, every homomorphism $\varphi \colon B \to X$ in $\V$ arises from a split extension of $B$ by~$X$ in $\V$. To this end, consider the split extension
\begin{equation}\label{split_phi}
\begin{tikzcd}
X \arrow[r, "\iota_2"] &
B \ltimes X \arrow[r, shift left, "\pi_1"] &
B, \arrow[l, shift left, "\iota_1"]
\end{tikzcd}
\end{equation}
where $B \ltimes X$ is the algebra whose underlying vector space is $B \times X$ and whose multiplication is defined by
\[
(b,x) \cdot (b',x') = (bb', xx' + \varphi(b)x' + x\varphi(b')),
\]
and $\iota_1, \iota_2$ and $\pi_1$ are the canonical injections and projection. Straightforward computations show that $\tau_B$ sends \eqref{split_phi} to the morphism
\[
\Inn \circ \, \varphi \colon B \to \E(X) \colon b \mapsto (L_{\varphi(b)}, R_{\varphi(b)}).
\]
Consequently, $\widetilde{\tau}_B$ sends the split extension \eqref{split_phi} to the morphism $\varphi$.

Finally, we need to check that \eqref{split_phi} is indeed a split extension in $\V$. The linear map
\[
\theta \colon B \ltimes X \to B \times X \colon (b,x) \mapsto (b, \varphi(b)+x)
\]
is bijective, and it is an algebra homomorphism onto the direct product $B \times X$, since
\begin{align*}
\theta(b,x) \cdot \theta(b',x') &= (bb', \varphi(b)\varphi(b') + \varphi(b)x' + x\varphi(b') + xx')\\
&= (bb', \varphi(bb') + xx' + \varphi(b)x' + x\varphi(b')) = \theta\big((b,x) \cdot (b',x')\big).
\end{align*}
Hence $B \ltimes X \cong B \times X$ is an algebra of $\V$.

Therefore, $\widetilde{\tau}_B$ is surjective, and the actions on $X$ are representable. Since $\Inn$ is an isomorphism, it follows that $X$ is its own actor.
\end{proof}

\begin{example}\label{ex_assoc}
Let $\V=\Assoc$ and fix the standard choice of the $\lambda/\mu$-rules. Then, for any associative algebra $X$, the map
\[
\Inn \colon X \to \Bim(X) \colon x \mapsto (L_x,R_x)
\]
is an isomorphism if and only if $X$ is a unital algebra.

Indeed, if $\Inn$ is an isomorphism, then $X$ is a unital algebra with unit
\[
1=\Inn^{-1}(\id_X,\id_X).
\]
Conversely, if $X$ has unit $1$, then $\Ker(\Inn)=\Ann(X)=\{ 0 \}$, and $\Inn$ is injective. Furthermore, as observed in~\cite[Remark~3.8]{WRAAlg}, any bimultiplier of $X$ is uniquely determined by the choice of an element $x=f \ast 1 = 1 \ast f \in X$, since
\[
f \ast y = xy \quad \text{and} \quad y \ast f = yx
\]
for any $y \in X$. Thus, $\Inn$ is surjective and, by \Cref{prop_e(x)}, there is a natural isomorphism
\[
\SplExt(-,X) \cong \Hom_\Assoc(-,X).
\]
Hence, actions on a unital associative algebra $X$ are representable. Similarly, one can show that the same result holds for the variety $\CAssoc$ of commutative associative algebras, replacing $\Bim(X)$ with the associative algebra $\M(X)$ of multipliers of $X$.
\end{example}

\begin{example}\label{ex_jjord}
Let $\V=\JJord$ be the variety of Jacobi--Jordan algebras over a field of characteristic $3$, and let $X$ be an object of $\V$. We choose the $\lambda/\mu$-rules
\[
x(yz)=-(xy)z-(xz)y, \qquad (yz)x=-(xy)z-(xz)y.
\]
We recall from~\cite{WRAAlg,SurveyActions} that, in this case, $\E(X)$ is an external actor of $X$, and it is isomorphic to the partial algebra $(\ADer(X),D_X,\mu_X)$, where
\[
\ADer(X)=\{ f \in \End(X) \mid f(xy) = -f(x)y - f(y)x, \; \forall x,y \in X \}
\]
is the vector space of \emph{anti-derivations} of $X$, the linear map $\widetilde m_X\colon\ADer(X)\otimes\ADer(X)\to\End(X)$ is determined by $\widetilde m_X(f\otimes g)=-f \circ g - g \circ f$, and
\[
D_X=\widetilde m_X^{-1}\bigl(\ADer(X)\bigr), \qquad \mu_X=\left.\widetilde m_X\right|_{D_X}.
\]
Thus, the product of $f,g\in\ADer(X)$ is defined if and only if the \emph{anti-commutator} $\{f,g\}=-f\circ g-g\circ f$ is an anti-derivation, and in this case it is equal to $\{f,g\}$. Again, we prove that the homomorphism of partial algebras
\[
\Inn \colon \widetilde{U}(X) \to \ADer(X) \colon x \mapsto L_x
\]
is an isomorphism if and only if $X$ is a unital Jacobi--Jordan algebra. Indeed, since $\bchar(\F)=3$, the identity map $\id_X$ is an anti-derivation of $X$, because $\id_X(xy)=xy=-2xy=-xy-yx$. Hence, if $\Inn$ is an isomorphism, then $X$ is a unital algebra with unit
\[
1=\Inn^{-1}(\id_X).
\]
Conversely, if $X$ has unit $1$, then $\Ker(\Inn)=\Ann(X)=\{ 0 \}$, and $\Inn$ is injective. Furthermore, any anti-derivation of $X$ is uniquely determined by the choice of an element $x=f(1)$, since
\[
f(y)=f(y \cdot 1)=-f(y)1-f(1)y
\]
for any $y \in X$. Hence $2f(y)=-f(1)y$ and, since $\bchar(\F)=3$, we get
\[
f(y)=f(1)y=xy.
\]
Thus, $\Inn$ is surjective. Since $\Inn$ is a bijective homomorphism of partial algebras with total domain, it is an isomorphism of partial algebras by \Cref{lem_iso_total}, and, by \Cref{prop_e(x)}, the actions on a unital Jacobi--Jordan algebra $X$ are representable.
\end{example}

In \Cref{SecMain}, we will show that the above proposition applies to all unital algebras in any operadic, action accessible, unit-closed variety $\V$.

\section{Representability of actions in ideally exact categories}\label{SecIdEx}

Let $\V$ be a unit-closed variety. In order to study the representability of actions of unital non-associative algebras, we need to recall the definition of an \emph{ideally exact category}, introduced by G.~Janelidze in~\cite{IdeallyExact}.
 
\begin{definition}\cite{IdeallyExact}
Let $\U$ be a category with pullbacks, with initial object $0$ and terminal object $1$. Then $\U$ is said to be \textit{ideally exact} if it is Barr-exact, Bourn-protomodular, has finite coproducts and the unique morphism $0 \to 1$ in $\U$ is a regular epimorphism.
\end{definition}

In addition to all semi-abelian categories, examples of ideally exact categories include the categories $\Ring$ and $\CRing$ of unital rings and unital commutative rings, all non-trivial protomodular varieties of universal algebras, such as the variety of \emph{MV-algebras}~\cite{rel}, and any cotopos. 

Moreover, the following characterisation holds.
 
\begin{theorem}\cite{IdeallyExact}\label{carid}
Let $\U$ be a category with pullbacks. The following conditions are equivalent:
\begin{enumerate}
\item $\U$ is ideally exact.
\item $\U$ is Barr-exact, has finite coproducts and there exists a monadic functor $\U \to \V$, where $\V$ is a semi-abelian category.
\item There exists a monadic functor $\U \to \V$, where $\V$ is a semi-abelian category, such that the underlying functor of the corresponding monad preserves regular epimorphisms and kernel pairs. \noproof
\end{enumerate}
\end{theorem}

\begin{remark}\label{cartesian}
Let $\U$ be an ideally exact category. It follows from \Cref{carid} that there exists a monadic adjunction 
\begin{equation}\label{monadic}
\begin{tikzcd}
{\U} & {\V}
\arrow[""{name=0, anchor=center, inner sep=0}, "U"', from=1-1, to=1-2]
\arrow[""{name=1, anchor=center, inner sep=0}, "F"', curve={height=12pt}, from=1-2, to=1-1]
\arrow["\dashv"{anchor=center, rotate=-90}, draw=none, from=1, to=0]
\end{tikzcd}
\end{equation}
with $\V$ semi-abelian. By~\cite[Theorem~3.3]{IdeallyExact}, this adjunction is associated with the unique morphism $0 \to 1$ (up to an equivalence) if and only if the unit of the adjunction is cartesian. 

Of course the monadic adjunction associated to $0 \to 1$ is always there. Indeed, one may take $\V = (\U \downarrow 0)$, the slice category over $0$, where $U=!^*$ is the pullback functor along $0 \to 1$, and $F$ is the functor which sends a morphism $X \to 0$ to its domain. However, this construction is not always the most natural choice. In particular, if $\U$ is already semi-abelian, it may be preferable to take $\V=\U$. We further observe that, since $F$ is a left adjoint, $F(0)$ is an initial object of $\U$.
\end{remark}

\begin{example}\label{monadic_ring}
Let $\U=\Ring$, which has the ring of integers $\Z$ as initial object and the zero ring $\{0\}$ as terminal one. Then a monadic adjunction as in \Cref{cartesian} is given by
\begin{equation*}\label{eq:situation_rng}
\begin{tikzcd}
{\Ring} & {\Rng,}
\arrow[""{name=0, anchor=center, inner sep=0}, "U"', from=1-1, to=1-2]
\arrow[""{name=1, anchor=center, inner sep=0}, "F"', curve={height=12pt}, from=1-2, to=1-1]
\arrow["\dashv"{anchor=center, rotate=-90}, draw=none, from=1, to=0]
\end{tikzcd}
\end{equation*}
where $U$ is the forgetful functor and $F$ maps any (not necessarily unital) ring $X$ to the semidirect product~$\Z \ltimes X$ with operations
\[
(n,x)+(m,y)=(n+m,x+y), \quad (n,x) \cdot (m,y)=(nm,xy+ny+mx)
\]
and unit element $(1,0)$. The unit $\eta \colon 1_{\Rng} \Rightarrow UF$ of the adjunction is cartesian since the map $\eta_X \colon X \to  U(\Z \ltimes X) \colon x \mapsto (0,x)$ is the kernel of
\[
U(\pi_1) \colon U(\Z \ltimes X) \to  U(\Z) \colon (n,x) \mapsto n.
\]
Thus, $F \dashv U$ is, up to an equivalence, the adjunction associated with the unique morphism $\Z \to \{ 0 \}$ in $\Ring$.

Similarly, there is a monadic adjunction with cartesian unit
\[
\begin{tikzcd}
{\CRing} & {\CRng.}
\arrow[""{name=0, anchor=center, inner sep=0}, "U"', from=1-1, to=1-2]
\arrow[""{name=1, anchor=center, inner sep=0}, "F"', curve={height=12pt}, from=1-2, to=1-1]
\arrow["\dashv"{anchor=center, rotate=-90}, draw=none, from=1, to=0]
\end{tikzcd}
\]
\end{example}

\begin{example}\label{monadic_alg}
Let $\V$ be a unit-closed variety. The class $\V_1$ of unital algebras of~$\V$ can be seen as a category whose morphisms are the algebra homomorphisms preserving the unit. Then $\V_1$ is a protomodular variety of universal algebras. Hence, by~\cite[Example~3.6]{IdeallyExact}, it is an ideally exact category. Moreover, $\V_1$ is not pointed, since an initial object is the field $\F$ endowed with the standard multiplication, while the terminal one is the zero algebra $\{0\}$.

Furthermore, there exists a monadic adjunction with cartesian unit
\[
\begin{tikzcd}
{\V_1} & {\V,}
\arrow[""{name=0, anchor=center, inner sep=0}, "U"', from=1-1, to=1-2]
\arrow[""{name=1, anchor=center, inner sep=0}, "F"', curve={height=12pt}, from=1-2, to=1-1]
\arrow["\dashv"{anchor=center, rotate=-90}, draw=none, from=1, to=0]
\end{tikzcd}
\]
where $U$ is the forgetful functor and $F$ maps any algebra $X$ of $\V$ to the semidirect product $\F \ltimes X$, with operations
\[
(n,x)+(m,y) = (n+m, x+y), \quad (n,x)\cdot(m,y) = (nm, xy + ny + mx)
\]
and multiplicative unit $(1,0)$. As in \Cref{monadic_ring}, the unit $\eta$ of the adjunction is cartesian, since $\eta_X \colon X \to U(\F \ltimes X) \colon x \mapsto (0,x)$ is the kernel of $U(\pi_1) \colon U(\F \ltimes X) \to U(\F) \colon (n,x) \mapsto n$.
\end{example}

In~\cite{IdeallyExact}, G.~Janelidze redefined the notion of action representability in the setting of ideally exact categories.

\begin{definition}\cite{IdeallyExact}\label{def_act_repr}
Let $\U$ be an ideally exact category and let $!^* \colon \U \to (\U \downarrow 0)$ be the pullback functor along $0 \to 1$. Then $\U$ is \emph{action representable} if the functor
\[
\SplExt(-,!^*(X)) \cong \Act(-,!^*(X)) \colon (\U \downarrow 0)^{\op} \to \Set
\]
is representable for any object $X$ of $\U$.
\end{definition}

\begin{example}\cite{IdeallyExact}\label{ex_ring}
Let $\V=\Rng$ be the (semi-abelian) category of not necessarily unital rings. As it is shown in~\cite{BJK2}, actions on a ring $X$ are representable in some cases, in particular when $X$ is unital. In this case, the canonical ring homomorphism
\[
X \to \Bim(X) \colon x \mapsto (L_x,R_x),
\]
where $\Bim(X)$ is the ring of \emph{bimultipliers} of $X$~\cite{MacLane58}, is an isomorphism. 

Now, recall there is a canonical category equivalence $(\Ring \downarrow \Z) \simeq \Rng$, which is defined as follows: for any object $X \to \Z$ of $(\Ring \downarrow \Z)$, its kernel is an object of~$\Rng$; conversely, every not necessarily unital ring $X$ can be seen as the kernel of the projection $\Z \ltimes X \to \Z \colon (n,x) \mapsto n$.

Under this equivalence, the pullback functor $!^* \colon \Ring \to (\Ring \downarrow \Z)$ along the unique morphism $\Z \to \{ 0 \}$ is the same as the forgetful functor $U \colon \Ring \to \Rng$. Indeed, given any unital ring $X$, viewed as the object $X \to \{ 0 \}$ of $(\Ring \downarrow \{ 0 \})$, its image under $!^*$ is given by the pullback
\[
\begin{tikzcd}
\Z \times X \arrow[r, "\pi_2"] \arrow[d, swap, "\pi_1"] \arrow[dr, phantom, "\usebox\pullback" , very near start, color=black] & X \arrow[d] \\
\Z \arrow[r] & \{ 0 \}
\end{tikzcd}
\]
and one has that $\Ker \pi_1 \cong X$, seen as a ring without unit. Since, under the equivalence $(\Ring \downarrow \Z) \simeq \Rng$, any object $A \to \Z$ corresponds to its kernel, it follows that $!^*$ acts by forgetting the unit.

Thus, given any unital ring $X$, the functor $\SplExt(-, !^{*}(X)) \colon (\Ring \downarrow \Z)^{\op} \to \Set$ is representable if and only if so is the functor $\SplExt(-,U(X)) \colon \Rng^{\op} \to  \Set$. Hence, as observed in~\cite[Section 5]{IdeallyExact}, we conclude that $\Ring$ is action representable, and that the actor of a unital ring $X$ is $X$ itself. The same argument applies to the category $\CRing$ of unital commutative rings.
\end{example}

\Cref{ex_ring} suggests that, in order to study the representability of actions of the category of unital rings, it is more convenient to choose the monadic adjunction with cartesian unit of \Cref{monadic_ring}.  In a similar way, for any unit-closed variety of non-associative $\F$-algebras $\V$, there is a category equivalence $(\V_1 \downarrow \F) \simeq \V$. Thus, actions on a unital algebra $X$ are representable if and only if the functor
\[
\SplExt(-,U(X)) \colon \V^{\op} \to \Set,
\]
where $U \colon \V_1 \to \V$ denotes the forgetful functor, is representable.

\begin{remark}
It is straightforward to verify that, if a unit-closed variety $\V$ is action representable, then so is the ideally exact category $\V_1$ of its unital algebras. This applies, for instance, to the category $\BoolAlg$ of \emph{Boolean algebras}, which can be described as the subcategory of Boolean rings with unit and unit-preserving morphisms; see~\cite{Bell}.

It was shown in~\cite[Proposition~3.1]{BJK2} that the variety of Boolean rings is action representable, and that the actor of a Boolean ring $X$ is the ring of $X$-linear endomorphisms of $X$. In particular, it follows that actions on a unital Boolean ring~$X$ are representable.

Moreover, in this case, the actor of a Boolean algebra $X$ can be identified with~$X$ itself. Indeed, any $X$-linear endomorphism $f \colon X \to X$ is uniquely determined by the element $x = f(1) \in X$, since for every $y \in X$ one has $f(y) = f(1)y = xy$.

Thus, the ideally exact category $\BoolAlg$ is action representable, and the actor of any Boolean algebra $X$ is $X$ itself.
\end{remark}

Following~\cite{IdeallyExact}, one may also re-define the notions of action accessibility and weak action representability.

\begin{definition}\label{act_acc}
Let $\U$ be an ideally exact category and let $!^* \colon \U \to (\U \downarrow 0)$ be the pullback functor along $0 \to 1$. Then
\begin{itemize}
\item $\U$ is \emph{action accessible} if, for any object $X$ of it, every object in the category $\SplExt(!^*(X))$ of split extensions with kernel $!^*(X)$ is accessible.
\item $\U$ is \emph{weakly action representable} if the functor $\SplExt(-,!^*(X))$ admits a weak representation for any object $X$ of $\U$.
\end{itemize}
\end{definition}

\begin{remark}
We observe that the chain of implications
\[
\textit{action representability}\Rightarrow\textit{weak action representability}\Rightarrow\textit{action accessibility}
\]
still holds in the setting of ideally exact categories. In fact:
\begin{itemize}
\item If $\U$ is action representable, then it is immediate that it is weakly action representable.
\item If $\U$ is weakly action representable, then Corollary~4.3 and Proposition~4.4 of~\cite{WRA} apply, in particular, to all split extensions in $(\U \downarrow 0)$ with kernel of the form $!^*(X)$, for some object $X$ of $\U$. Consequently, every object in $\SplExt(!^*(X))$ is accessible, and hence $\U$ is an action accessible category. We point the reader to~\cite[Theorem~2.22]{SurveyActions}, where an explicit proof of this implication for semi-abelian categories is given.
\end{itemize}
\end{remark}

In the next section, we prove that, in the context of unital non-associative algebras, the existence of the $\lambda/\mu$-rules implies action representability. Our main tool is to show that the existence of a multiplicative unit is equivalent to the canonical map $\Inn \colon \widetilde{U}(X) \to \E(X)$ being an isomorphism.

Before doing so, we show that action accessibility fails in general in the setting of categories of unital algebras, as the following proposition illustrates. In this setting, as in the case of action representability, we replace the pullback functor along $\F \to \{0\}$ with the forgetful functor $U \colon \V_1 \to \V$. In other words, we work with the monadic adjunction with cartesian unit described in \Cref{monadic_alg}.

\begin{proposition}
The ideally exact categories $\Alg_1$ and $\Com_1$ are not action accessible.
\end{proposition}

\begin{proof}
Let $\V=\Com$ be the category of commutative (not necessarily associative) algebras, and consider the monadic adjunction of \Cref{monadic_alg}. We construct a commutative algebra $B$ and a split extension of $B$ by $U(\F)$ in $\Com$, where $\F$ is the initial object of $\Com_1$, which does not admit any morphism into a faithful object of $\SplExt(U(\F))$. Our construction is inspired by the argument in~\cite[Proposition~6.4]{Cappelletti}, where it is shown that the category of \emph{lattice-ordered groups} is not action accessible. The idea is to adapt the technique used there to the setting of non-associative algebras.

We take $B$ to be the $6$-dimensional commutative algebra with basis
\[
\{ 1,e_1,e_2,e_3,f_1,f_2 \}, 
\]
where $1$ is a multiplicative unit of $B$, and whose multiplication table is defined by
\[
e_i^2=f_j^2=1, \quad e_if_j=f_je_i=0, \quad f_1 f_2 = f_2 f_1=e_1
\]
and
\[
e_1 e_2 = e_2 e_1=e_3, \quad e_2e_3=e_3e_2=e_1, \quad e_3e_1=e_1e_3=e_2.
\]
We observe that $B$ is a simple algebra, i.e., its only ideals are $\{ 0 \}$ and $B$ itself. Indeed, let $I$ be a non-zero ideal of $B$, and let
\[
0 \neq b=a1+\sum_{i=1}^3 b_i e_i + \sum_{j=1}^2 c_j f_j \in I,
\]
with $a,b_1,b_2,b_3,c_1,c_2 \in \F$. Then
\[
((bf_1)e_1)f_1=((af_1+c_1 1+c_2e_1)e_1)f_1=(c_1e_1+c_2 1)f_1=c_2f_1 \in I.
\]
If $c_2 \neq 0$, then $1=f_1^2 \in I$ and $I=B$. If instead $c_2=0$, then
\[
(b f_1)e_1=c_1 e_1 \in I.
\]
Again, if $c_1 \neq 0$, then $1=e_1^2 \in I$ and $I=B$. If instead $c_1=c_2=0$, then
\[
b f_1=a f_1 \in I.
\]
If $a \neq 0$, then $1 \in I$ and $I=B$. Now, suppose that $a=c_1=c_2=0$. In this case
\[
b=b_1 e_1 + b_2 e_2 + b_3 e_3
\]
with $b_i \neq 0$ for some $i \in \{1,2,3\}$. Since $e_je_i \in \operatorname{span}\{e_1,e_2,e_3\}$ for $j \neq i$ and $e_kf_1=0$ for every $k$, we get
\[
((be_i)f_1)f_1=(b_if_1)f_1=b_i 1 \in I,
\]
which implies that $I=B$. Thus, $B$ is a simple algebra.

Now let $\sigma \colon B \to B$ be the linear map exchanging $f_1$ and $f_2$ and fixing $1,e_1,e_2,e_3$. It is immediate from the multiplication table that $\sigma$ is an algebra automorphism of~$B$, and $\sigma \neq \id_B$. Let moreover $\varepsilon \colon B \to \F$ be the linear form sending each element of $B$ to its coefficient along $1$ in the basis $\{ 1,e_1,e_2,e_3,f_1,f_2 \}$. Then $\varepsilon(1)=1$ and $\varepsilon \circ \sigma = \varepsilon$.

Consider the action of $B$ on $U(\F)$ given by
\[
b \ast t = t \ast b = \varepsilon(b)t,
\]
for any $b \in B$ and $t \in \F$. Since $\Com$ is determined by commutativity alone, the algebra $B \ltimes_\varepsilon U(\F)$, whose underlying vector space is $B \times \F$ and whose multiplication is defined by
\[
(b,t) \cdot (b',t') = (bb', tt' + \varepsilon(b)t' + \varepsilon(b')t),
\]
is commutative. Hence we obtain a split extension
\begin{equation}\label{non_faithful}
\begin{tikzcd}
U(\F) \arrow [r, "\iota_2"]
& B \ltimes_\varepsilon U(\F) \arrow[r, shift left, "\pi_1"] &
B \ar[l, shift left, "\iota_1"]
\end{tikzcd}
\end{equation}
in $\Com$, where $\iota_1,\iota_2$ and $\pi_1$ are the canonical injections and projection. Since $\varepsilon \circ \sigma = \varepsilon$, the linear map $\sigma \times \id_\F$ is an algebra endomorphism of $B \ltimes_\varepsilon U(\F)$, and we get two distinct morphisms in $\SplExt(U(\F))$:
\[
\begin{tikzcd}
U(\F) \arrow [r, "\iota_2"] \ar[d, equal]
& B \ltimes_\varepsilon U(\F) \arrow[r, shift left, "\pi_1"] \ar[d, "\id_B \times \id_\F"']&
B \ar[l, shift left, "\iota_1"]\ar[d, "\id_B"] \\
U(\F) \arrow [r, "\iota_2"]
& B \ltimes_\varepsilon U(\F) \arrow[r, shift left, "\pi_1"] &
B \ar[l, shift left, "\iota_1"]
\end{tikzcd}
\]
and
\[
\begin{tikzcd}
U(\F) \arrow [r, "\iota_2"] \ar[d, equal]
& B \ltimes_\varepsilon U(\F) \arrow[r, shift left, "\pi_1"] \ar[d, "\sigma \times \id_\F"']&
B \ar[l, shift left, "\iota_1"]\ar[d, "\sigma"] \\
U(\F) \arrow [r, "\iota_2"]
& B \ltimes_\varepsilon U(\F) \arrow[r, shift left, "\pi_1"] &
B \ar[l, shift left, "\iota_1"]
\end{tikzcd}
\]
We deduce that \eqref{non_faithful} is not faithful.

So, if $\Com_1$ were action accessible, then there should exist a faithful object
\[
\begin{tikzcd}
U(\F) \arrow [r, "k'"]
& A' \arrow[r, shift left, "\alpha'"] &
B' \ar[l, shift left, "\beta'"]
\end{tikzcd}
\]
in $\SplExt(U(\F))$ and a morphism
\[
\begin{tikzcd}
U(\F) \arrow [r, "\iota_2"] \ar[d, equal]
& B \ltimes_\varepsilon U(\F) \arrow[r, shift left, "\pi_1"] \ar[d, "g"']&
B \ar[l, shift left, "\iota_1"]\ar[d, "f"] \\
U(\F) \arrow [r, "k'"]
& A' \arrow[r, shift left, "\alpha'"] &
B' \ar[l, shift left, "\beta'"] .
\end{tikzcd}
\]
Since $\Ker f$ is an ideal of $B$ and $B$ is simple, we have two possibilities: $f$ is injective or $f$ is the zero morphism.

If $f$ is injective, then \eqref{non_faithful} is faithful (see~\cite[Proposition~1.4]{act_accessible}). Indeed, if $(g_1,f_1)$ and $(g_2,f_2)$ are two morphisms in $\SplExt(U(\F))$ with codomain \eqref{non_faithful} and the same domain, then $(g \circ g_1, f \circ f_1)=(g \circ g_2, f \circ f_2)$, since the codomain of $(g,f)$ is faithful. Hence $f_1=f_2$, because $f$ is injective, and $g_1=g_2$, since the middle component of a morphism in $\SplExt(U(\F))$ is determined by its components on the kernel and on the base. This contradicts the fact that \eqref{non_faithful} is not faithful.

If $f$ is the zero morphism, then, for any $b \in B$ and $t \in \F$, we have
\[
k'(\varepsilon(b)t) = g(\iota_2(\varepsilon(b)t)) = g\big((b,0) \cdot (0,t)\big) = g(\iota_1(b)) \cdot g(\iota_2(t)) = \beta'(f(b)) \cdot k'(t) = 0.
\]
Since $k'$ is injective, it follows that $\varepsilon(b)t=0$ for any $b \in B$ and $t \in \F$, which contradicts $\varepsilon(1)=1$.

Thus, we conclude that $\Com_1$ is not action accessible. Furthermore, the split extension \eqref{non_faithful} is also a split extension in $\Alg$, and the same argument applies verbatim. Hence the category $\Alg_1$ of all unital not necessarily associative algebras is not action accessible either.
\end{proof}

\section{Action representability of unital algebras}\label{SecMain}

The aim of this section is to prove that, given an operadic, action accessible, unit-closed variety of non-associative algebras $\V$, the ideally exact category $\V_1$ is action representable. We start by proving the following proposition which establishes a connection between action accessibility and the property of being a unit-closed variety.

\begin{proposition}\label{lemma_unit}
Let $\V$ be an operadic, action accessible, unit-closed variety of non-associative algebras, and fix a choice of the $\lambda/\mu$-rules.
\begin{enumerate}
\item If $\V$ satisfies a non-trivial identity of degree $2$, then
\begin{itemize}
\item if $\bchar(\F)\neq 3$, $\V=\CAssoc$;
\item if $\bchar(\F)= 3$, $\V=\CAssoc$ or $\V$ is a subvariety of $\JJord$.
\end{itemize}
\item If $\V$ satisfies no identities of degree $2$, then 
\begin{equation}\label{eq_lemma1}
\begin{alignedat}{2}
\lambda_1 &= 1 + \lambda_4, \qquad & \lambda_2 &= -\lambda_4 - \lambda_7 - \lambda_8, \\
\lambda_3 &= -\lambda_4 - \lambda_6 - \lambda_8, \qquad & \lambda_5 &= \lambda_8
\end{alignedat}
\end{equation}
and
\begin{equation}\label{eq_lemma2}
\begin{alignedat}{2}
\mu_1 &= \mu_4, \qquad & \mu_2 &= 1 - \mu_4 - \mu_7 - \mu_8, \\
\mu_3 &= 1 - \mu_4 - \mu_6 - \mu_8, \qquad & \mu_5 &= -1 + \mu_8.
\end{alignedat}
\end{equation}
\item For any algebra $X$ of $\V$, the pair $(\id_X,\id_X)$ belongs to $\E(X)$ and it is the unit of the partial algebra $\E(X)$.
\end{enumerate}
\end{proposition}

\begin{proof}{\ }
\begin{enumerate}
\item Suppose that $\V$ satisfies a non-trivial identity of degree $2$. By~\cite[Proposition~1.3]{Tim}, $\V$ is a subvariety of $\Com$ or of $\ACom$. Since $\V$ is unit-closed, the field $\F$, endowed with its multiplication, is an object of $\V$. If $\V$ is a subvariety of $\ACom$, then $1\cdot1=-1\cdot1$, so that $\bchar(\F)=2$ and $\ACom=\Com$. Hence, in both cases, $\V$ is a subvariety of $\Com$. Consequently, the $\lambda/\mu$-rules reduce to a single equation of the form
\[
x(yz)=\alpha(xy)z+\beta(xz)y
\]
for some $\alpha,\beta \in \F$. By~\cite[Proposition~2.1]{WRAAlg}, it follows that $\V$ is either a subvariety of $\CAssoc$, or a subvariety of $\JJord$. If $\bchar(\F)\neq 3$, the second case is excluded, since the field $\F$, endowed with its multiplication, is an object of every unit-closed variety but is not a Jacobi--Jordan algebra. Moreover, the only proper operadic subvarieties of $\CAssoc$ are the varieties $\Nil_k(\CAssoc)$, for $k\geq1$, and the trivial variety, none of which is unit-closed. Therefore, if $\V$ is a subvariety of $\CAssoc$, then $\V = \CAssoc$.  Note that, if $\bchar(\F)=3$, then $\CAssoc$ is itself a subvariety of $\JJord$, since $x(yz)+y(zx)+z(xy)=3xyz=0$ in every commutative associative algebra.

\item Suppose that $\V$ does not satisfy any identity of degree $2$, and let $X$ be the free algebra on two generators $a,b$. Since $\V$ is unit-closed, the unital algebra $\widetilde X=\langle X,1\rangle$ is again an object of $\V$. Applying the $\lambda/\mu$-rules in $\widetilde X$ with $(x,y,z)=(1,a,b)$, we obtain
\[
ab=(\lambda_1+\lambda_2+\lambda_7+\lambda_8)\,ab +(\lambda_3+\lambda_4+\lambda_5+\lambda_6)\,ba
\]
and
\[
ab=(\mu_1+\mu_2+\mu_7+\mu_8)\,ab +(\mu_3+\mu_4+\mu_5+\mu_6)\,ba.
\]
Since there are no identities of degree $2$, then $ab$ and $ba$
are linearly independent, and it follows that
\begin{equation}\label{system1}
\begin{alignedat}{2}
\lambda_1+\lambda_2+\lambda_7+\lambda_8 &= 1, \qquad &
\lambda_3+\lambda_4+\lambda_5+\lambda_6 &= 0, \\
\mu_1+\mu_2+\mu_7+\mu_8 &= 1, \qquad &
\mu_3+\mu_4+\mu_5+\mu_6 &= 0.
\end{alignedat}
\end{equation}
Similarly, by applying the $\lambda/\mu$-rules respectively with
$(x,y,z)=(a,1,b)$ and $(x,y,z)=(a,b,1)$, one gets
\begin{equation}\label{system2}
\begin{alignedat}{2}
\lambda_1+\lambda_2+\lambda_5+\lambda_7 &= 1, \qquad &
\lambda_3+\lambda_4+\lambda_6+\lambda_8 &= 0, \\
\mu_3+\mu_4+\mu_6+\mu_8 &= 1, \qquad &
\mu_1+\mu_2+\mu_5+\mu_7 &= 0
\end{alignedat}
\end{equation}
and
\begin{equation}\label{system3}
\begin{alignedat}{2}
\lambda_1+\lambda_3+\lambda_5+\lambda_6 &= 1, \qquad &
\lambda_2+\lambda_4+\lambda_7+\lambda_8 &= 0, \\
\mu_2+\mu_4+\mu_7+\mu_8 &= 1, \qquad &
\mu_1+\mu_3+\mu_5+\mu_6 &= 0.
\end{alignedat}
\end{equation}
Solving the system given by equations
\eqref{system1}--\eqref{system2}--\eqref{system3}, one obtains
\eqref{eq_lemma1} and \eqref{eq_lemma2}.

\item Let $X$ be a not necessarily unital algebra of $\V$. We first claim that $(\id_X,\id_X)\in\E(X)$. If $X$ is unital, then
\[
(\id_X,\id_X) = (L_{1_X},R_{1_X}) = \Inn(1_X) \in \E(X).
\]
If $X$ does not have a multiplicative unit, we consider the unital algebra $\widetilde X=\langle X,1 \rangle$, which is an object of $\V$ since $\V$ is unit-closed. Then
\[
(\id_{\widetilde X},\id_{\widetilde X}) = (L_1,R_1) = \Inn(1) \in \E(\widetilde X).
\]
By definition of $\E(\widetilde X)$, this means that for every identity $\Phi(x_1,\ldots,x_k)=0$ which determines $\V$ we have
\begin{equation}\label{eq_identity}
\Phi(\alpha_1,\ldots,\alpha_k)=0,
\end{equation}
for any choice of $\alpha_j=\id_{\widetilde X}$ and $\alpha_t\in\widetilde X$, with $t \neq j \in \{1,\ldots,k\}$. Since $X \subseteq \widetilde X$, equation~\eqref{eq_identity} holds in particular for all $\alpha_t\in X$. Therefore $(\id_X,\id_X)$ satisfies all the identities of $\V$, that is, $(\id_X,\id_X) \in \E(X)$.

Next, we observe that, since $\V$ is unit-closed, evaluating the $\lambda/\mu$-rules in $\widetilde X$ at the triples $(x,y,1)$ and $(x,1,y)$, with $x,y\in X$, shows that the multilinear identities
\begin{align}
xy&=(\lambda_1+\lambda_3+\lambda_5+\lambda_6)\,xy+(\lambda_2+\lambda_4+\lambda_7+\lambda_8)\,yx, \label{deg2_1}\\
xy&=(\lambda_1+\lambda_2+\lambda_5+\lambda_7)\,xy+(\lambda_3+\lambda_4+\lambda_6+\lambda_8)\,yx, \label{deg2_2}\\
yx&=(\mu_1+\mu_3+\mu_5+\mu_6)\,xy+(\mu_2+\mu_4+\mu_7+\mu_8)\,yx, \label{deg2_3}\\
yx&=(\mu_1+\mu_2+\mu_5+\mu_7)\,xy+(\mu_3+\mu_4+\mu_6+\mu_8)\,yx \label{deg2_4}
\end{align}
hold in every algebra of $\V$. By \Cref{rem_one_f}, they still hold when $x\in X$ and $y$ is replaced by an element $f\in\E(X)$.

Now let $f=(f\ast-,-\ast f)\in\E(X)$, and set $e=(\id_X,\id_X)$. By definition of $m_X$, for every $x\in X$ we have
\begin{align*}
x\ast m_X(f,e) &= (\lambda_1+\lambda_3+\lambda_5+\lambda_6)(x\ast f)+(\lambda_2+\lambda_4+\lambda_7+\lambda_8)(f\ast x),\\
m_X(f,e)\ast x &= (\mu_1+\mu_3+\mu_5+\mu_6)(x\ast f)+(\mu_2+\mu_4+\mu_7+\mu_8)(f\ast x),\\
x\ast m_X(e,f) &= (\lambda_1+\lambda_2+\lambda_5+\lambda_7)(x\ast f)+(\lambda_3+\lambda_4+\lambda_6+\lambda_8)(f\ast x),\\
m_X(e,f)\ast x &= (\mu_1+\mu_2+\mu_5+\mu_7)(x\ast f)+(\mu_3+\mu_4+\mu_6+\mu_8)(f\ast x).
\end{align*}
By \eqref{deg2_1}--\eqref{deg2_4} with $y=f$, the right-hand sides are equal to $x\ast f$, $f\ast x$, $x\ast f$ and $f\ast x$, respectively. Hence
\[
m_X(f,e)=m_X(e,f)=f\in\E(X),
\]
so that $f\otimes e,\,e\otimes f\in D_X$ and $\langle f,e\rangle=\langle e,f\rangle=f$. Thus, $e$ is the unit of the partial algebra $\E(X)$.
\end{enumerate}
\end{proof}

\begin{example}
Let $\V=\Assoc$ and fix the standard choice of the $\lambda/\mu$-rules as in \Cref{Example associative}. Then, one may check that \eqref{eq_lemma1}--\eqref{eq_lemma2} are satisfied.
\end{example}

\begin{example}\label{example Alt}
Let $\bchar(\F) \neq 2$ and consider the variety $\Alt$ of alternative algebras. The multilinearisation process shows that $\Alt$ is equivalent to the variety determined by
\[
(xy)z+(xz)y-x(yz)-x(zy)=0
\]
and
\[
(xy)z+(yx)z-x(yz)-y(xz)=0.
\] 
Conversely, setting $z=y$ in the first identity and $y=x$ in the second one gives twice the right and the left alternative identities, respectively. This is where the assumption $\bchar(\F)\neq2$ is used. A choice of the $\lambda/\mu$-rules is then given by
\[
\lambda_1=\lambda_2=-\lambda_7=-\mu_2=\mu_7=\mu_8=1,
\]
and all remaining coefficients equal to zero. Thus, \eqref{eq_lemma1} and \eqref{eq_lemma2} are satisfied.
\end{example}

We observe that the results of \Cref{lemma_unit} do not hold for a variety $\V$ that is not unit-closed. For instance, if $\V = \Lie$, then $\E(X) \cong \Der(X)$. However, the identity map $\id_X$ is a derivation of $X$ if and only if $X$ is abelian. 

Moreover, the conditions \eqref{eq_lemma1}--\eqref{eq_lemma2} are not sufficient to ensure that a variety of non-associative algebras $\V$ satisfying no identities of degree~$2$ is unit-closed. For instance, let $\V = \Nil_3(\Assoc)$ and fix the choice of the $\lambda/\mu$-rules as in \Cref{Example associative}. Then \eqref{eq_lemma1}--\eqref{eq_lemma2} are satisfied, but $\V$ is not unit-closed. Nevertheless, the following characterisation holds.

\begin{corollary}\label{cor_unit}
Let $\V$ be an operadic, action accessible variety of non-associative algebras, and fix a choice of the $\lambda/\mu$-rules. Suppose moreover that $\V$ satisfies no identities of degree $2$ and that all the identities of $\V$ are consequences of the $\lambda/\mu$-rules. Then $\V$ is unit-closed if and only if \eqref{eq_lemma1}--\eqref{eq_lemma2} hold.
\end{corollary}

\begin{proof}
If $\V$ is unit-closed, then \eqref{eq_lemma1}--\eqref{eq_lemma2} hold by \Cref{lemma_unit}.

Conversely, assume that \eqref{eq_lemma1}--\eqref{eq_lemma2} hold. We must show that $\V$ is unit-closed, i.e., for any algebra $X$ in $\V$, the algebra $\widetilde X=\langle X,1 \rangle$ is still an object of~$\V$.

Since every identity of $\V$ is a consequence of the $\lambda/\mu$-rules, it is sufficient to verify that \eqref{eq00}--\eqref{eq01} hold for all $x,y,z \in \widetilde X$. These identities already hold when $x,y,z \in X$, since $X$ is an algebra of $\V$. Thus, it remains to check them in the cases where at least one of the variables is the external unit $1$. That is, it is enough to consider triples of the form
\[
(1,y,z),\ (x,1,z),\ (x,y,1),\ (1,1,z),\ (1,y,1),\ (x,1,1),\ (1,1,1).
\]
In each of these cases, the identities \eqref{eq00}--\eqref{eq01} reduce to trivial equalities as a consequence of \eqref{eq_lemma1}--\eqref{eq_lemma2}. This completes the proof.
\end{proof}

We are now ready to prove our main result, namely that, given an operadic, unit-closed variety $\V$, the existence of $\lambda/\mu$-rules is sufficient to guarantee that the ideally exact category $\V_1$ is action representable. 

We begin by observing that, for a unit-closed variety $\V$, the existence of the $\lambda/\mu$-rules in $\V_1$ is enough to ensure the existence of the $\lambda/\mu$-rules in~$\V$. Indeed, assume that there exist $\lambda_1,\ldots,\lambda_8,\mu_1,\ldots,\mu_8 \in \F$ such that the $\lambda/\mu$-rules hold in $\V_1$. Then, any (not necessarily unital) algebra $X$ of $\V$ can be identified with the kernel of the projection
\[
\pi_1 \colon \F \ltimes X \to \F \colon (n,x) \mapsto n.
\]
Hence, $X$ is (isomorphic to) an ideal of the unital algebra $\F \ltimes X$, and therefore the $\lambda/\mu$-rules hold on $X$ as well. It follows that $\V$ satisfies the $\lambda/\mu$-rules and, by \Cref{Theorem AC iff Orzech}, we have the following.

\begin{proposition}
Let $\V$ be an operadic, unit-closed variety of non-associative algebras. The following are equivalent:
\begin{tfae}
\item $\V$ is action accessible;
\item The $\lambda/\mu$-rules hold in $\V$;
\item The $\lambda/\mu$-rules hold in $\V_1$. \noproof
\end{tfae}
\end{proposition}

We observe that these conditions are not equivalent to the action accessibility of the ideally exact category~$\V_1$; see \Cref{prop_jord}.

\begin{theorem}\label{thm_unital}
Let $\V$ be an operadic, action accessible, unit-closed variety of non-associative algebras, and fix a choice of the $\lambda/\mu$-rules. Let $X$ be an object of $\V$ and let $\widetilde{U} \colon \V \to \PAlg$ be the forgetful functor. The map
\[
\Inn \colon \widetilde{U}(X) \to \E(X) \colon x \mapsto (L_x,R_x)
\]
is an isomorphism of partial algebras if and only if $X$ is a unital algebra.
\end{theorem}

\begin{proof}
If $\Inn$ is an isomorphism of partial algebras, then $X$ is a unital algebra with unit
\[
1=\Inn^{-1}(\id_X,\id_X),
\]
since $(\id_X,\id_X)$ is the unit of $\E(X)$ by \Cref{lemma_unit}(3).

Conversely, suppose that $X$ has a multiplicative unit $1$. Then $\Ker(\Inn)=\Ann(X)=\{ 0 \}$, and $\Inn$ is injective. To show that $\Inn$ is surjective, let $f=(f \ast -, - \ast f) \in \E(X)$ and $x \in X$. By \Cref{rem_one_f}, we may evaluate the $\lambda/\mu$-rules at triples containing $f$ once and elements of $X$ otherwise, where as before $fx=f \ast x$ and $xf=x \ast f$. Applying \eqref{eq00} with $(x,y,z)=(x,f,1)$, we get
\[
x(f \ast 1)=(\lambda_1+\lambda_3+\lambda_5+\lambda_6)(x \ast f)+(\lambda_2+\lambda_4+\lambda_7+\lambda_8)(f \ast x),
\]
and the right-hand side is equal to $x \ast f$ by \eqref{deg2_1} with $y=f$. Thus
\[
x \ast f=x(f \ast 1)
\]
for every $x \in X$, and for $x=1$ we obtain $f \ast 1=1 \ast f$. Similarly, applying \eqref{eq01} with $(x,y,z)=(x,1,f)$ and using \eqref{deg2_4} with $y=f$, we get
\[
f \ast x=(1 \ast f)x.
\]
This implies that
\[
f=(L_a,R_a)=\Inn(a),
\]
where $a=f \ast 1 = 1 \ast f \in X$, and $\Inn$ is surjective. Since $\Inn$ is a bijective homomorphism of partial algebras with total domain, it is an isomorphism of partial algebras by \Cref{lem_iso_total}.
\end{proof}

As a consequence, we get the following.

\begin{theorem}\label{final_thm}
Let $\V$ be an operadic, action accessible, unit-closed variety of non-associative algebras. The ideally exact category $\V_1$ of unital algebras of $\V$ is action representable, and the actor of a unital algebra $X$ is $X$ itself.
\end{theorem}

\begin{proof}
It follows from \Cref{prop_e(x)} and \Cref{thm_unital} that, for any unital algebra $X$ of $\V$, there exists a natural isomorphism
\[
\SplExt(-,U(X)) \cong \Hom_\V(-,U(X)),
\]
where $U \colon \V_1 \rightarrow \V$ is the forgetful functor. Thus, $\V_1$ is action representable, and $X$ is its own actor.
\end{proof}

\begin{examples}{\ }
\begin{enumerate}
\item The ideally exact categories $\Assoc_1$ and $\CAssoc_1$ of unital associative and unital commutative associative algebras are action representable. Indeed, as shown in \Cref{ex_assoc}, for any unital associative algebra $X$ (resp.~unital commutative associative algebra), the map $\Inn$ induces an isomorphism $\Bim(X) \cong X$ (resp.~$\M(X) \cong X$).

\item
Let $\bchar(\F) \neq 2$ and consider the variety $\Alt$ of alternative algebras. By \Cref{example Alt}, $\Alt$ is operadic and action accessible, and it is unit-closed. Hence, by \Cref{final_thm}, the ideally exact category $\Alt_1$ of unital alternative algebras is action representable, that is, for any unital alternative algebra~$X$, there is a natural isomorphism
\[
\SplExt(-,U(X)) \cong \Hom_\Alt(-,U(X)),
\]
as also observed in~\cite[Example~3.11(8)]{WRAAlg}.

\item Let $\WAlt$ be the variety determined by the identity
\begin{equation}\label{xabi}
x(yz)-(xy)z - z(xy) + (zx)y = 0.
\end{equation}
We call $\WAlt$ the variety of \emph{weak alternative algebras}. One may check that a choice of the $\lambda/\mu$-rules is given by
\[
\lambda_1=\lambda_3=-\lambda_6=-\mu_3=\mu_6=\mu_8=1,
\]
and hence $\WAlt$ is an operadic, action accessible variety. Moreover, $\WAlt$ satisfies no identities of degree $2$, its defining identity \eqref{xabi} is one of the $\lambda/\mu$-rules, and \eqref{eq_lemma1}--\eqref{eq_lemma2} hold. Thus, $\WAlt$ is unit-closed by \Cref{cor_unit}.

Now, the identity \eqref{xabi} can be rewritten as
\[
[x,y,z]-[z,x,y]=0,
\]
where $[x,y,z]=x(yz)-(xy)z$ denotes the \emph{associator}. It is well known that the associator of any alternative algebra $X$ is totally skew-symmetric, that is
\[
[x_{\sigma(1)},x_{\sigma(2)},x_{\sigma(3)}]=\operatorname{sgn}(\sigma)[x_1,x_2,x_3],
\]
for any permutation $\sigma \in S_3$. In particular, since the cycle $(1\,3\,2)$ is even, we obtain
\[
[z,x,y]=[x,y,z].
\]
Therefore, every alternative algebra satisfies \eqref{xabi}, and hence $\Alt \subseteq \WAlt$. Note that this inclusion still holds in characteristic $2$.

Now, by \Cref{final_thm}, the ideally exact category $\WAlt_1$ is action representable, and the actor of any algebra $X$ of $\WAlt_1$ is $X$ itself. Consequently, the category $\Alt_1$ of unital alternative algebras is action representable also in characteristic $2$. Indeed, let $X$ be a unital alternative algebra and $B$ an alternative algebra. By \Cref{prop_e(x)} and \Cref{thm_unital}, every split extension of $B$ by $U(X)$ in $\WAlt$ is isomorphic to $B \ltimes U(X)$, constructed as in~\eqref{split_phi} from a unique homomorphism $\varphi \colon B \to U(X)$, and $B \ltimes U(X) \cong B \times U(X)$ is an alternative algebra (see the proof of \Cref{prop_e(x)}). Hence $\SplExt(-,U(X)) \cong \Hom_{\Alt}(-,U(X))$.

\item Let $\bchar(\F)=3$ and consider the variety $\JJord$. Then, the ideally exact category $\JJord_1$ of unital Jacobi--Jordan algebras is action representable. In fact, as observed in \Cref{ex_jjord}, the partial algebra $\ADer(X)$ of anti-derivations of $X$ is isomorphic to $X$ precisely when $X$ is unital.
\end{enumerate}
\end{examples}

The result of \Cref{final_thm} highlights a contrast with the non-unital case. Indeed, in the setting of semi-abelian varieties of non-associative algebras, as shown in~\cite{Tim}, action representability is a highly restrictive property, essentially characterising the variety of Lie algebras. Here, actions on $X$ become representable whenever $X$ is a unital algebra of an operadic, action accessible, unit-closed variety.

The action accessibility of $\V$ is sufficient, but not necessary, for the action representability of~$\V_1$, as the following result shows. To the best of our knowledge, it provides the first example of an action representable ideally exact category whose associated semi-abelian category, through the monadic adjunction of \Cref{monadic_alg}, is not even action accessible.

\begin{proposition}\label{prop_jord}
Let $\bchar(\F)=0$. Then the variety $\Jord$ of Jordan algebras is unit-closed and it is not action accessible. Nevertheless, the ideally exact category $\Jord_1$ is action representable, and the actor of a unital Jordan algebra $X$ is $X$ itself.
\end{proposition}

\begin{proof}
We have shown in \Cref{SecVar} that $\Jord$ is unit-closed. Now, since the Jordan identity has degree~$4$, every multilinear identity of degree~$3$ of $\Jord$ is a consequence of commutativity. Hence, the monomials $x(yz)$, $y(xz)$ and $z(xy)$ are linearly independent in the free Jordan algebra on $x,y,z$, while, by commutativity, the right-hand side of~\eqref{eq00} is a linear combination of $y(xz)$ and $z(xy)$. Thus, $\Jord$ satisfies no $\lambda/\mu$-rules and, since it is homogeneous, it is not action accessible by \Cref{Theorem AC iff Orzech}. Alternatively, the variety $\Jord$ is operadic, since $\bchar(\F)=0$ (see \Cref{SecVar}), and unit-closed. Hence, if it were action accessible, then \Cref{lemma_unit}(1) would imply $\Jord=\CAssoc$. This is false: the Jordan algebra of symmetric $2\times2$ matrices, with product $a\circ b=(ab+ba)/2$, is not associative. For instance, if $p=E_{11}$ and $s=E_{12}+E_{21}$, then $(p\circ s)\circ s=I_2/2$, whereas $p\circ(s\circ s)=p$.

We now prove that $\Jord_1$ is action representable. Since $\bchar(\F)=0$, the multilinearisation process (see~\cite[Corollary~3.7]{Osborn}) shows that every Jordan algebra satisfies the \emph{linearised Jordan identity}
\begin{equation}\label{eq_lin_jord}
(ab)(cy)+(bc)(ay)+(ca)(by)=c((ab)y)+a((bc)y)+b((ca)y),
\end{equation}
that is, $[L_a,L_{bc}]+[L_b,L_{ca}]+[L_c,L_{ab}]=0$. Let $X$ be a unital Jordan algebra with unit $1$, let $B$ be a Jordan algebra, and consider a split extension
\begin{equation}\label{split_jord}
\begin{tikzcd}
U(X) \arrow[r, "k"] & A \arrow[r, shift left, "\alpha"] & B \arrow[l, shift left, "\beta"]
\end{tikzcd}
\end{equation}
in $\Jord$. Set $e=k(1)$ and $A_0=\{a\in A \mid ea=0\}$. We claim that
\[
A=k(X)\oplus A_0,\qquad A_0A_0\subseteq A_0, \qquad k(X)A_0=0.
\]
Since $k(X)$ is an ideal of $A$ with unit $e$, for every $a\in A$ we have $ea\in k(X)$ and $e(ea)=ea$, so that $a-ea\in A_0$. Moreover, if $a\in k(X)\cap A_0$, then $a=ea=0$. This proves the first assertion. If $u,v\in A_0$, then evaluating~\eqref{eq_lin_jord} at $(a,b,c,y)=(e,e,u,v)$ gives $e(uv)=0$, which proves the second one. Finally, let $z\in k(X)$ and $u\in A_0$. Then $w=zu$ belongs to $k(X)$, so that $we=w$ and $e(we)=w$, and evaluating~\eqref{eq_lin_jord} at $(a,b,c,y)=(e,z,u,e)$ gives $w=2w$, that is, $zu=0$. This proves the claim.

Thus, $A$ is the direct product of its subalgebras $k(X)$ and $A_0$, and the projection
\[
\pi\colon A\to k(X)\colon a\mapsto ea
\]
is an algebra homomorphism. We define
\[
\Phi_B\colon \SplExt(B,U(X))\to\Hom_{\Jord}(B,U(X))
\]
by sending the class of~\eqref{split_jord} to the homomorphism $\varphi=k^{-1}\circ\pi\circ\beta$, that is, $\varphi(b)=k^{-1}(e\beta(b))$. This does not depend on the representative: if $\gamma$ is an isomorphism from~\eqref{split_jord} to a split extension
\[
\begin{tikzcd}
U(X) \arrow[r, "\bar k"] & \bar A \arrow[r, shift left, "\bar\alpha"] & B, \arrow[l, shift left, "\bar\beta"]
\end{tikzcd}
\]
and $\bar e=\bar k(1)$ and $\bar\pi$ are defined as above, then $\gamma(e)=\bar e$, so that $\gamma\circ\pi=\bar\pi\circ\gamma$, and hence
\[
\bar k^{-1}\circ\bar\pi\circ\bar\beta=\bar k^{-1}\circ\bar\pi\circ\gamma\circ\beta=\bar k^{-1}\circ\gamma\circ\pi\circ\beta=k^{-1}\circ\pi\circ\beta.
\]

Conversely, for a homomorphism $\varphi\colon B\to U(X)$, let $\Psi_B(\varphi)$ be the class of the split extension
\[
\begin{tikzcd}
U(X) \arrow[r, "\iota_2"] & B\times U(X) \arrow[r, shift left, "\pi_1"] & B, \arrow[l, shift left, "\beta_\varphi"]
\end{tikzcd}
\]
where $B\times U(X)$ is the direct product of algebras, $\iota_2$ and $\pi_1$ are the canonical injection and projection, and $\beta_\varphi(b)=(b,\varphi(b))$, which is a homomorphism since $\varphi$ is. In $B\times U(X)$ we have $e=(0,1)$ and $e(b,\varphi(b))=(0,\varphi(b))$, so that $\Phi_B(\Psi_B(\varphi))=\varphi$.

Now, let $\varphi$ be the image under $\Phi_B$ of the class of~\eqref{split_jord}, and let $\psi=\beta-k\circ\varphi$, that is, $\psi(b)=\beta(b)-e\beta(b)$. Since $\psi$ is the composite of $\beta$ with the projection $A\to A_0$, it is an algebra homomorphism with values in $A_0$. Moreover, $\alpha\circ\psi=\id_B$ and, for every $u\in A_0$, we have $\psi(\alpha(u))-u\in A_0\cap\Ker\alpha=A_0\cap k(X)=0$. Hence, $\psi\colon B\to A_0$ is an isomorphism, and so is the map
\[
\theta\colon B\times U(X)\to A\colon (b,x)\mapsto \psi(b)+k(x).
\]
Since $\theta\circ\iota_2=k$, $\alpha\circ\theta=\pi_1$ and $\theta\circ\beta_\varphi=\psi+k\circ\varphi=\beta$, the map $\theta$ is an isomorphism of split extensions, and $\Psi_B(\varphi)$ is the class of~\eqref{split_jord}. Therefore, $\Phi_B$ is a bijection with inverse~$\Psi_B$.

Finally, $\Phi_B$ is natural in $B$. Indeed, let $g\colon B'\to B$ be a morphism in $\Jord$. The pullback of~\eqref{split_jord} along $g$ is the split extension
\begin{equation}\label{split_jord2}
\begin{tikzcd}
U(X) \arrow[r, "k'"] & A' \arrow[r, shift left, "\alpha'"] & B', \arrow[l, shift left, "\beta'"]
\end{tikzcd}
\end{equation}
where
\[
A'=\{(a,b')\in A\times B' \mid \alpha(a)=g(b')\},
\]
with $k'(x)=(k(x),0)$, $\alpha'(a,b')=b'$ and $\beta'(b')=(\beta(g(b')),b')$. Then $e'=(e,0)$ and
\[
e'\beta'(b')=(e\beta(g(b')),0)=k'(\varphi(g(b'))),
\]
so that $\Phi_{B'}$ sends the class of \eqref{split_jord2} to $\varphi\circ g$. We conclude that
\[
\SplExt(-,U(X))\cong\Hom_{\Jord}(-,U(X)),
\]
that is, $\Jord_1$ is action representable, and $X$ is its own actor.
\end{proof}

\begin{remark}\label{rem_direct_factor}
In all the cases considered above, if $X$ is a unital algebra, then every split extension of an algebra $B$ by $U(X)$ is isomorphic to the direct product $B\times U(X)$, with section $b\mapsto(b,\varphi(b))$ for a unique homomorphism $\varphi\colon B\to U(X)$. Equivalently, $U(X)$ is a direct factor of the middle object of every split extension with kernel $U(X)$. This does not trivialise the split extensions: different homomorphisms $\varphi$ give non-isomorphic split extensions, and the corresponding action of $B$ on $X$ is given by $b\ast x=\varphi(b)x$ and $x\ast b=x\varphi(b)$. Conversely, the second part of the proof of \Cref{prop_jord} shows that, for any unit-closed variety $\V$, if $U(X)$ is a direct factor of the middle object of every split extension with kernel $U(X)$, then $X$ is its own actor. In this sense, unital algebras behave like complete groups, which are direct factors of every group containing them as normal subgroups. The split extension~\eqref{non_faithful} shows that this property fails in $\Com$ and in $\Alg$: if it were of this form, then $\varepsilon(b)=b\ast1=\varphi(b)$ would be an algebra homomorphism, whereas $\varepsilon(e_1^2)=1\neq0=\varepsilon(e_1)^2$.

The construction of the external weak actor provides a sufficient condition for this property, which is not necessary, as \Cref{prop_jord} shows. It remains an open problem to characterise the unit-closed varieties $\V$ for which the ideally exact category $\V_1$ is action representable, for instance by means of identities playing the role of the $\lambda/\mu$-rules. We do not know whether action representability, weak action representability and action accessibility coincide for the ideally exact categories~$\V_1$, nor, in particular, whether the action accessibility of $\V_1$ implies that, for every unital algebra $X$, the kernel $U(X)$ is a direct factor of the middle object of every split extension with kernel $U(X)$.
\end{remark}

\section{Action representability of unital Poisson algebras}\label{SecPois}

The aim of this section is to study the representability of actions of unital Poisson algebras. We prove that the categories $\Pois_1$ and $\CPois_1$ of unital Poisson algebras and unital commutative Poisson algebras over a field $\F$ with $\bchar(\F) \neq 2$ are action representable, by employing the explicit description of split extensions and derived actions of Poisson algebras given in~\cite{CigoliManciniMetere}, and the construction of the external weak actor given in~\cite[Section~5.2]{SeveralMult}. The assumption $\bchar(\F)\neq2$ is needed in order to apply these results, which are stated under this hypothesis. This also suggests a possible approach to extending the present analysis to ideally exact categories of unital algebras with several bilinear, not necessarily associative, multiplications, using the external weak actor constructed in~\cite[Section~4]{SeveralMult}.

\begin{definition}
A Poisson algebra is a vector space $X$ over a field $\F$ equipped with two bilinear operations
\[
\cdot \colon X \times X\to X \quad \text{and} \quad [-,-] \colon X \times X  \to X,
\]
such that $(X,\cdot)$ is an associative algebra, $(X,[-,-])$ is a Lie algebra and the \emph{Poisson identity} holds:
\[
[x,yz]=[x,y]z+y[x,z], \quad \forall x,y,z \in X.
\]
This means that for any $x \in X$, the adjoint map $\operatorname{ad}_x=[x,-]$ is a derivation of the associative algebra $(X,\cdot)$.

A Poisson algebra is said to be \emph{commutative} (resp.~\emph{unital}) if the underlying associative algebra is commutative (resp.~unital).
\end{definition}

\begin{remark}
Let $X$ be a Poisson algebra. We denote by
\[
\Ann(X)=\{ x \in X \mid xy=yx=0, \; \forall y \in X \}
\]
the annihilator of the associative algebra $(X,\cdot)$, and by
\[
\bZ(X)=\{ x \in X \mid [x,y]=0, \; \forall y \in X \}
\]
the centre of the Lie algebra $(X,[-,-])$.
\end{remark}

We observe that, if $X$ is a unital Poisson algebra, then the unit element $1$ belongs to the centre $\bZ(X)$. Indeed, applying the Poisson identity with $y=z=1$, one gets
\[
[x,1]=[x,1]1+1[x,1]=2[x,1]
\]
and, hence, $[x,1]=[1,x]=0$ for any $x \in X$.

\begin{example}
Any associative algebra $(X,\cdot)$ becomes a Poisson algebra with the usual Lie bracket
\[
[x,y]=xy-yx.
\]
Furthermore:
\begin{itemize}
\item If $(X,\cdot)$ is a commutative associative algebra, then $(X,\cdot,[-,-])$ is a commutative Poisson algebra with trivial Lie bracket.
\item If $(X,\cdot)$ is a unital associative algebra, then $(X,\cdot,[-,-])$ is a unital Poisson algebra.
\end{itemize}
\end{example}

We denote by $\Pois$ and $\CPois$ the semi-abelian varieties of (commutative) Poisson algebras, and by $\Pois_1$ and $\CPois_1$ the subcategories of unital (commutative) Poisson algebras. We observe that $\Pois$ and $\CPois$ are Orzech categories of interest, so they are action accessible by~\cite[Theorem~4.7]{Montoli}. However, they are not action representable, as shown in~\cite[Theorem~5.3]{SeveralMult}.

\begin{remark}
The categories $\Pois_1$ and $\CPois_1$ are ideally exact, and their initial object is the field $\F$ endowed with the standard multiplication and the trivial Lie bracket $[n,m]=0$. Furthermore, there is a monadic adjunction
\[
\begin{tikzcd}
{\Pois_1} & {\Pois,}
\arrow[""{name=0, anchor=center, inner sep=0}, "U"', from=1-1, to=1-2]
\arrow[""{name=1, anchor=center, inner sep=0}, "F"', curve={height=12pt}, from=1-2, to=1-1]
\arrow["\dashv"{anchor=center, rotate=-90}, draw=none, from=1, to=0]
\end{tikzcd}
\]
where $U$ is the forgetful functor and $F$ maps every not necessarily unital Poisson algebra $X$ to the semidirect product $\F \ltimes X$ with operations
\begin{align*}
(n,x)+(m,y)&=(n+m,x+y),\\
(n,x)\cdot(m,y)&=(nm, xy + ny + mx),\\
[(n,x),(m,y)]&=(0,[x,y])
\end{align*}
and unit $(1,0)$. As in the case of associative algebras, the unit $\eta \colon 1_{\Pois} \Rightarrow UF$ of the adjunction is cartesian since the map $\eta_X \colon X \to  U(\F \ltimes X) \colon x \mapsto (0,x)$ is the kernel of
\[
U(\pi_1) \colon U(\F \ltimes X) \to  U(\F) \colon (n,x) \mapsto n.
\]
Thus, $F \dashv U$ is, up to an equivalence, the adjunction associated with the unique morphism $\F \to \{ 0 \}$ in $\Pois_1$.

Similarly, there is a monadic adjunction with cartesian unit
\[
\begin{tikzcd}
{\CPois_1} & {\CPois.}
\arrow[""{name=0, anchor=center, inner sep=0}, "U"', from=1-1, to=1-2]
\arrow[""{name=1, anchor=center, inner sep=0}, "F"', curve={height=12pt}, from=1-2, to=1-1]
\arrow["\dashv"{anchor=center, rotate=-90}, draw=none, from=1, to=0]
\end{tikzcd}
\]
\end{remark}

Let $X$ be a Poisson algebra. Following~\cite[Section~5.2]{SeveralMult}, we denote by $\E(X)$ the vector space of all triples $f=(f \ast -,- \ast f,[f,-]) \in \Bim(X)\times \Der(X)$ such that
\begin{align*}
f \ast [x,y]&=[f \ast x,y]-[f,y]x,\\
[x,y] \ast f&=[x \ast f,y]-x[f,y],\\
[f,xy]&=[f,x]y+x[f,y],
\end{align*}
for all $x,y \in X$, endowed with the bilinear operations
\[
f \cdot g=(f \ast (g \ast -), (- \ast f) \ast g, f \ast [g,-] + [f,-] \ast g)
\]
and
\[
[f,g]=(f \ast [g,-] - [g,f \ast -], [g,-] \ast f - [g,- \ast f], [f,[g,-]]-[g,[f,-]]).
\]
In~\cite{SeveralMult}, these operations are denoted by $\langle f,g\rangle_\ast$ and $\langle f,g\rangle_{[-,-]}$, respectively. By~\cite[Theorem~5.6]{CigoliManciniMetere}, $\E(X)$ is an object of the category $\Alg^2$ of algebras with two not necessarily associative bilinear operations, and there exists a natural monomorphism of functors
\[
\tau \colon \SplExt(-,X) \rightarrowtail \Hom_{\Alg^2}(\overline{U}(-),\E(X)),
\]
where $\overline{U} \colon \Pois \to \Alg^2$ is the forgetful functor. By~\cite[Section~5.2]{SeveralMult}, $\E(X)$ is the external weak actor of~$X$, and $\tau$ is the corresponding external weak representation.

Moreover, a morphism 
\[
\varphi \colon \overline{U}(B) \to \E(X) \colon b \mapsto (b \ast -, - \ast b, [b,-])
\]
in $\Alg^2$ belongs to $\Imm(\tau_B)$ if and only if 
\begin{equation}\label{eq_pois}
(b \ast x) \ast b' = b \ast (x \ast b'),
\end{equation}
for any $b,b' \in B$ and $x \in X$ (see~\cite[Theorem~5.6]{CigoliManciniMetere}).

\begin{remark}
Let $X$ be a Poisson algebra. We observe that the non-associative algebra $(\E(X),\cdot)$ has multiplicative unit
\[
1=(\id_X,\id_X,0_X),
\]
where $0_X \colon X \to X$ denotes the zero map. Indeed, one has
\[
f \cdot 1=1 \cdot f=f,
\]
for any $f \in \E(X)$. Furthermore, $1 \in \bZ(\E(X))$, since $[f,1]=[1,f]=0$.
\end{remark}

We further recall from~\cite[Remark~5.8]{CigoliManciniMetere} that, for any Poisson algebra $X$, if the underlying associative algebra $(X,\cdot)$ has trivial annihilator or is perfect, that is, $(X^2,\cdot)=(X,\cdot)$, then
\[
f \ast (x \ast g) = (f \ast x) \ast g,
\]
for every $f=(f \ast -, - \ast f)$ and $g=(g \ast -,- \ast g)$ in $\Bim(X)$ and $x \in X$.

Consequently, for any other Poisson algebra $B$, every morphism from $\overline{U}(B)$ to $\E(X)$ lies in $\Imm(\tau_B)$, and hence there is a natural isomorphism
\begin{equation}\label{isopois}
\SplExt(-,X)\cong \Hom_{\Alg^{2}}(\overline{U}(-),\E(X)).
\end{equation}

Since any unital associative algebra is perfect and with trivial annihilator, the isomorphism \eqref{isopois} holds in particular for any object $X$ of $\Pois_1$. Moreover, we can prove the following.

\begin{proposition}\label{prop_pois}
Let $X$ be a Poisson algebra. The algebra homomorphism
\[
\Inn \colon \overline{U}(X) \to \E(X) \colon x \mapsto (L_x, R_x, \operatorname{ad}_x),
\] 
is an isomorphism if and only if $X$ is a unital algebra.
\end{proposition}

\begin{proof}
If $\Inn$ is an isomorphism, then $X$ is a unital Poisson algebra with
\[
1=\Inn^{-1}(\id_X,\id_X,0_X).
\]
Conversely, suppose that $X$ is a unital Poisson algebra. We have
\[
\Ker(\Inn)=\Ann(X)\cap\bZ(X)=\{0\},
\]
and $\Inn$ is injective. It remains to show that $\Inn$ is surjective. To this end, let
\[
f=(f \ast -, - \ast f, [f,-]) \in \E(X).
\]
Since $(X,\cdot)$ is unital, by \Cref{ex_assoc} we have an isomorphism of associative algebras $\Bim(X) \cong X$, and there exists $a \in X$ such that
\[
f \ast x= ax \quad \text{and} \quad x \ast f= xa
\]
for any $x \in X$. Moreover, since $a=f \ast 1=1 \ast f$, one has
\[
0=f \ast [1,x]=[f \ast 1,x]-[f,x]1=[a,x]-[f,x],
\]
and hence $[f,x]=[a,x]$ for every $x \in X$. Thus, $f=(L_a,R_a,\operatorname{ad}_a)=\Inn(a)$, and we conclude that $\Inn$ is an isomorphism.
\end{proof}

It follows that for any unital Poisson algebra $X$, there is a natural isomorphism
\[
\SplExt(-, U(X)) \cong \Hom_\Pois(-,U(X)),
\]
where $U \colon \Pois_{1} \to \Pois$ denotes the forgetful functor, and we get the following.

\begin{theorem}
Let $\bchar(\F)\neq 2$. The ideally exact category $\Pois_{1}$ is action representable, and the actor of a unital Poisson algebra $X$ is $X$ itself. \noproof
\end{theorem}

We observe that the same result holds in the category $\CPois_{1}$ of unital commutative Poisson algebras. In this case, by~\cite[Section~5.2]{SeveralMult}, the external weak actor of $X$ in $\CPois$ is isomorphic to the algebra
\[
\left\{
(f \ast -, [f,-])\in\M(X)\times\Der(X)
\ \middle|\
\begin{aligned}
f \ast [x,y]&=[f \ast x,y]-[f,y]\cdot x,\\
[f,x \cdot y]&=[f,x] \cdot y+x \cdot [f,y],
\end{aligned}
\; \; \forall x,y\in X
\right\},
\]
with two bilinear multiplications
\[
f\cdot g=(f \ast (g \ast -), f \ast [g,-] + [f,-] \ast g)
\]
and
\[
[f,g]=(f \ast [g,-] - [g,f \ast -], [f,[g,-]]-[g,[f,-]]).
\]
Here, the external weak representation is a natural isomorphism, so that this algebra is an external actor of~$X$; see~\cite[Section~5]{CigoliManciniMetere}. Moreover, the argument of \Cref{prop_pois} applies verbatim: if $X$ is unital, then $\M(X) \cong X$ via $f \mapsto f\ast 1$ (see \Cref{ex_assoc}), and the condition $f\ast[1,x]=[f\ast 1,x]-[f,x]$ yields $[f,-]=\operatorname{ad}_{f\ast 1}$. Hence, we conclude with the following.

\begin{theorem}
Let $\bchar(\F)\neq 2$. The ideally exact category $\CPois_{1}$ is action representable, and the actor of a unital commutative Poisson algebra $X$ is $X$ itself. \noproof
\end{theorem}

\section*{Acknowledgements}
The first and second authors are supported by the ``National Group for Algebraic and Geometric Structures and their Applications'' (GNSAGA -- INdAM). The first author was also supported by the SDF Sustainability Decision Framework Research Project -- MISE decree of 31/12/2021 (MIMIT Dipartimento per le politiche per le imprese -- Direzione generale per gli incentivi alle imprese) -- CUP:~B79J23000530005, COR:~14019279, Lead Partner:~TD Group Italia Srl, Partner:~University of Palermo, and he is a Postdoctoral Researcher of the Fonds de la Recherche Scientifique--FNRS. The third author was supported by Wallonie-Bruxelles International.


\section*{Use of generative AI tools}
The generative AI tool Anthropic's Claude Opus~5.5 was used to polish the exposition and to help check the correctness of the proofs; it was not used to obtain any mathematical results. The authors reviewed all suggestions made by this tool and take full responsibility for the content of this article.

\end{document}